\documentclass[12pt]{article}
\usepackage{amsmath,amssymb,amsfonts}
\usepackage{graphicx}
\usepackage{epsfig}
\usepackage{amsfonts}
\usepackage{amssymb, latexsym, amsmath, pb-diagram}
\usepackage{pstricks-add}
\usepackage{multicol}
\usepackage{bm}
\usepackage{mathrsfs}
\usepackage{xy}
\usepackage{epic,eepic}
\xyoption{all}
\begin{document}
\title{The Complement of Polyhedral Product Spaces and the Dual Simplicial Complexes}
\author {
Qibing Zheng\\School of Mathematical Science and LPMC, Nankai University\\
Tianjin 300071, China\\
zhengqb@nankai.edu.cn\footnote{Project Supported by Natural Science
Foundation of China, grant No. 11071125 and No. 11671154\newline
\hspace*{5.5mm}Key words and phrases: complement space; dual complex;
polyhedral join complex; universal algebra.\newline
\hspace*{5.5mm}Mathematics subject classification: 55N10}}\maketitle
\input amssym.def
\newsymbol\leqslant 1336
\newsymbol\geqslant 133E
\baselineskip=20pt
\def\w{\widetilde}
\def\RR{\mathscr R}
\def\SS{\mathscr S}
\def\XX{\mathscr X}

\begin{abstract} In this paper, we define and prove basic properties
of complement polyhedral product spaces, dual complexes and polyhedral join complexes.
Then we compute the universal algebra of polyhedral join complexes under certain split conditions
and the Alexander duality isomorphism on certain polyhedral product spaces.
\end{abstract}

\hspace*{40mm}${\displaystyle{\bf Table\,\,of\,\,Contents}}$

Section 1\, Introduction

Section 2\, Complement Spaces, Dual Complexes and Polyhedral Join Complexes

Section 3\, Homology and Cohomology Group

Section 4\, Universal Algebra

Section 5\, Duality Isomorphism
\vspace{3mm}

\newtheorem{Definition}{Definition}[section]
\newtheorem{Theorem}{Theorem}[section]
\newtheorem{Lemma}{Lemma}[section]
\newtheorem{Example}{Example}[section]
\font\hua=eusm10 scaled\magstephalf
\def\MM{\mathscr M}
\def\LL{\mathscr L}
\def\zz{\Bbb Z}

\section{Introduction}\vspace{3mm}

The polyhedral product theory, especially the homotopy type of polyhedral product
spaces, is developing rapidly nowadays. The first known polyhedral product space was the moment-angle
complex introduced by Buchstaber and Panov \cite{BP} and
was widely studied by mathematicians in the area of toric topology and geometry
(see \cite{A},\cite{B},\cite{GM},\cite{GR},\cite{GJ},\cite{GS}).
Later on, the homotopy types of polyhedral product spaces were studied by
Grbi\'{c} and Theriault \cite {GR},\cite{GJ},\cite{GS}, Beben and Grbi\'{c} \cite{GB},
Bahri, Bendersky, Cohen and Gitler \cite{B1},\cite{B2},\cite{B3} and many others
(\cite {BS},\cite{DJ},\cite{DE}).
The cohomology ring of homology split polyhedral product spaces and the
cohomology algebra over a field of polyhedral product spaces were computed in \cite{Z}.

For a polyhedral product space ${\cal Z}(K;\underline{X},\underline{A})$,
is the complement space
$M^c=(X_1{\times}{\cdots}{\times}X_m)\setminus {\cal Z}(K;\underline{X},\underline{A})$
a polyhedral product space?
In Theorem~2.4, we show that $M^c={\cal Z}(K^\circ;\underline{X},\underline{A}^c)$, where $K^\circ$
is the dual complex of $K$ relative to $[m]$ and $A_k^c=X_k{\setminus}A_k$.
The moment-angle complex case of this theorem is obtained by Gruji\'{c} and Welker \cite{GW}.

Let ${\cal Z}(K;\underline{Y},\underline{B})$,
$(\underline{Y},\underline{B})=\{(Y_k,B_k)\}_{k=1}^m$ be the polyhedral product space
defined as follows. For each $k$, $(Y_k,B_k)$ is a pair of polyhedral product spaces given by ($s_k=n_1{+}{\cdots}{+}n_k$)
$$Y_k={\cal Z}(X_k;\underline{U_k},\underline{C_k}),\,\, B_k={\cal Z}(A_k;\underline{U_k},\underline{C_k}),\,\,
(\underline{U_k},\underline{C_k})=\{(U_i,C_i)\}_{i=s_{k-1}{+}1}^{s_k},$$
where $(X_k,A_k)$ is a simplicial pair on $[n_k]$.
In Theorem 2.9, we prove that
${\cal Z}(K;\underline{Y},\underline{B})$ is also a polyhedral product space
${\cal Z}({\cal S}(K;\underline{X},\underline{A});\underline{U},\underline{C})$,
where ${\cal S}(K;\underline{X},\underline{A})$ is defined in Definition 2.7.
This theorem is a trivial extension of Ayzenberg's Proposition 5.1 of \cite{AY}, where
${\cal S}(K;\underline{X},\underline{A})$ is denoted by ${\cal Z}^*_K(\underline{X},\underline{A})$
and given the name polyhedral join complex.
When all $(X_k,A_k)=(\Delta\!^{[n_k]},K_k)$, the complex ${\cal S}(K;\underline{X},\underline{A})$ is
denoted by ${\cal S}(K;K_1,{\cdots},K_m)$ which is just
the composition complex $K(K_1,{\cdots},K_m)$ in Definition 4.5 of \cite{AY}.

In Section 3,
we compute the (co)homology group of polyhedral product inclusion complex $C_*(K;\underline{\vartheta})$ (defined in Definition 3.6)
in Theorem 3.7.
As an application, the reduced (co)homology group and the (right) total (co)homology group of polyhedral join complexes
is computed in Theorem~3.9 and Theorem~3.11.
In Example~3.10, we show that ${\cal S}(K;L_1,{\cdots},L_m)$ is a homological sphere
if and only if $K$ is a homological sphere and each $L_k=\partial\Delta\!^{[n_k]}$.
This result is a part of Ayzenberg's Theorem~6.6 of \cite{AY}, where homological can be replaced by simplicial.

In section~4, we compute the cohomology algebra of polyhedral product inclusion complex $C_*(K;\underline{\vartheta})$ in Theorem 4.7.
As an application, we compute the cohomology algebra of ${\cal Z}(K;\underline{Y},\underline{B})$ mentioned above
and the (right) universal (normal, etc.) algebra of ${\cal S}(K;\underline{X},\underline{A})$
(${\cal S}(K;L_1,{\cdots},L_m)$) in Example~4.8 (4.9).
In Theorem 5.6, we compute the Alexander duality isomorphism on the pair $
(X_1{\times}{\cdots}{\times}X_m,{\cal Z}(K;\underline{X},\underline{A}))$,
where all $X_k$'s are orientable manifolds and all $A_k$'s are polyhedra.
\vspace{3mm}

\section{Complement Spaces, Dual Complexes
and Polyhedral Join Complexes}\vspace{3mm}

{\bf Conventions and Notations} For a finite set $S$, $\Delta\!^S$ is the simplicial complex with only one maximal simplex $S$,
i.e., it is the set of all subsets of $S$ including the empty set $\emptyset$.
Define $\partial \Delta\!^S=\Delta\!^S{\setminus}\{S\}$.
For $[m]=\{1,{\cdots},m\}$, $\partial\Delta\!^{[m]}=\Delta\!^{[m]}{\setminus}[m]$.
Specifically, define
$\Delta\!^\emptyset=\{\emptyset\}$ and $\partial\Delta\!^{\emptyset}=\{\,\}$.
The void complex $\{\,\}$ with no simplex is inevitable in this paper.

For a simplicial complex $K$ on $[m]$ (ghost vertex $\{i\}\notin K$ is allowed) and $\sigma\subset[m]$ ($\sigma\notin K$ is allowed),
the link of $\sigma$ with respect to $K$ is the simplicial complex
${\rm link}_K\sigma=\{\tau\,|\,\sigma{\cup}\tau\in K,\,\sigma{\cap}\tau=\emptyset\}$.
This implies ${\rm link}_K\sigma=\{\emptyset\}$ if $\sigma$ is a maximal simplex of $K$ and ${\rm link}_K\sigma=\{\,\}$ if $\sigma\notin K$.
Specifically, if $K=\{\,\}$, then ${\rm  link}_K\sigma=\{\,\}$ for all $\sigma$.
\vspace{3mm}

{\bf Definition 2.1}
For a simplicial complex $K$ on $[m]$
and a sequence of topological (not CW-complex!) pairs $(\underline{X},\underline{A})=\{(X_k,A_k)\}_{k=1}^m$,
the {\it polyhedral product space} ${\cal Z}(K;\underline{X},\underline{A})$
is the subspace of $X_1{\times}{\cdots}{\times}X_m$ defined as follows.
For a subset $\tau$ of $[m]$, define
$$D(\tau)= Y_1{\times}\cdots{\times}Y_m,\quad Y_k=\left\{\begin{array}{cl}
X_k&{\rm if}\,\,k\in \tau, \\
A_k&{\rm if}\,\,k\not\in \tau.
\end{array}
\right.$$
Then ${\cal Z}(K;\underline{X},\underline{A})=\cup_{\tau\in K}\,D(\tau)$.
Empty space $\emptyset$ is allowed in a topological pair and $\emptyset{\times}X=\emptyset$ for all $X$.
Define ${\cal Z}(\{\,\};\underline{X},\underline{A})=\emptyset$.
\vspace{3mm}

Notice that $D(\sigma)={\cal Z}(\Delta\!^\sigma;\underline{X},\underline{A})$,
$D(\emptyset)=A_1{\times}{\cdots}{\times}A_m={\cal Z}(\{\emptyset\};\underline{X},\underline{A})$ and
$D([m])=X_1{\times}{\cdots}{\times}X_m={\cal Z}(\Delta\!^{[m]};\underline{X},\underline{A})$.
But $\emptyset={\cal Z}(\{\,\};\underline{X},\underline{A})$ has no corresponding $D(-)$.
\vspace{3mm}

{\bf Example 2.2} For ${\cal Z}(K;\underline{X},\underline{A})$,
let $S=\{k\,|\,A_k=\emptyset\}$. Then
$${\cal Z}(K;\underline{X},\underline{A})={\cal Z}({\rm link}_KS;\underline{X'},\underline{A'}){\times}(\Pi_{k\in S}X_k),$$
where $(\underline{X'},\underline{A'})=\{(X_k,A_k)\}_{k\notin S}$ and link is as defined in conventions.
\vspace{3mm}

{\bf Definition 2.3} Let $K$ be a simplicial complex with vertex set a subset of $S\neq\emptyset$.
The {\it dual of $K$ relative to} $S$ is the simplicial
complex
$$K^\circ=\{\,S{\setminus}\sigma\,\,|\,\,\sigma\subset S,\,\sigma\!\notin\! K\,\}.$$

It is obvious that
$(K^\circ)^\circ=K$, $(K_1{\cup}K_2)^\circ=(K_1)^\circ{\cap}(K_2)^\circ$ and
$(K_1{\cap}K_2)^\circ=(K_1)^\circ{\cup}(K_2)^\circ$. Specifically, $(\Delta\!^S)^\circ=\{\,\}$
and $(\partial\Delta\!^S)^\circ=\{\emptyset\}$.
\vspace{3mm}

{\bf Theorem 2.4} {\it For ${\cal Z}(K;\underline{X},\underline{A})$, the complement space
$${\cal Z}(K;\underline{X},\underline{A})^c
=(X_1{\times}{\cdots}{\times}X_m){\setminus}{\cal Z}(K;\underline{X},\underline{A})
={\cal Z}(K^\circ;\underline{X},\underline{A}^c),$$
where $(\underline{X},\underline{A}^c)=\{(X_k,A_k^c)\}_{k=1}^m$ with $A_k^c=X_k{\setminus}A_k$
and $K^\circ$ is the dual of $K$ relative to $[m]$.
}\vspace{2mm}

{\it Proof}\, For $\sigma\subset[m]$ but $\sigma\neq[m]$ ($\sigma=\emptyset$ is allowed),
$$\begin{array}{l}
\quad(X_1{\times}{\cdots}{\times}X_m)\setminus D(\sigma)\,\,({\rm with\,space\,pair}\,(X_k,A_k))\vspace{2mm}\\
=\cup_{j\notin\sigma}\,X_1{\times}{\cdots}{\times}(X_j{\setminus}A_j){\times}{\cdots}{\times}X_m\vspace{2mm}\\
=\cup_{j\in[m]\setminus\sigma}\,D([m]{\setminus}\{j\})\,\,({\rm with\,space\,pair}\,(X_k,A_k^c))\vspace{2mm}\\
={\cal Z}((\Delta\!^\sigma)^\circ;\underline{X},\underline{A}^c)
\end{array}$$
So for $K\neq\Delta\!^{[m]}$ or $\{\,\}$,
$$\begin{array}{l}
\quad{\cal Z}(\Delta\!^{[m]};\underline{X},\underline{A})\setminus{\cal Z}(K;\underline{X},\underline{A})\vspace{2mm}\\
={\cal Z}(\Delta\!^{[m]};\underline{X},\underline{A})\setminus\big(\cup_{\sigma\in K}{\cal Z}(\Delta\!^\sigma;\underline{X},\underline{A})\big)\vspace{2mm}\\
=\cap_{\sigma\in K}\big({\cal Z}(\Delta\!^{[m]};\underline{X},\underline{A})\setminus{\cal Z}(\Delta\!^\sigma;\underline{X},\underline{A})\big)\vspace{2mm}\\
=\cap_{\sigma\in K}{\cal Z}((\Delta\!^\sigma)^\circ;\underline{X},\underline{A}^c)\vspace{2mm}\\
={\cal Z}(K^\circ;\underline{X},\underline{A}^c)\end{array}$$

For $K=\Delta\!^{[m]}$ or $\{\,\}$, the above equality holds naturally. \hfill $\Box$
\vspace{3mm}

{\bf Example 2.5} Let $\Bbb F$ be a field and $V$ be a linear space over $\Bbb F$ with base $e_1,{\cdots},e_m$.
For a subset $\sigma=\{i_1,{\cdots},i_s\}\subset[m]$, denote by $\Bbb F(\sigma)$ the subspace of $V$ with base $e_{i_1},{\cdots},e_{i_s}$.
Then for $\Bbb F=\Bbb R$ or $\Bbb C$ and a simplicial complex $K$ on $[m]$, we have
$$V\setminus(\cup_{\sigma\in K}\Bbb R(\sigma))=\Bbb R^m\setminus{\cal Z}(K;\Bbb R,\{0\})={\cal Z}(K^\circ;\Bbb R,\Bbb R{\setminus}\{0\})
\simeq {\cal Z}(K^\circ;D^1,S^0),$$
$$V\setminus(\cup_{\sigma\in K}\Bbb C(\sigma))=\Bbb C^m\setminus{\cal Z}(K;\Bbb C,\{0\})={\cal Z}(K^\circ;\Bbb C,\Bbb C{\setminus}\{0\})
\simeq {\cal Z}(K^\circ;D^2,S^1).$$
This example is in Lemma 2.4 of \cite{GW}.
\vspace{3mm}

{\bf Theorem 2.6} {\it Let $K$ and $K^\circ$ be the dual of each other relative to $[m]$.
The index set $\XX_m=\{(\sigma\!,\,\omega)\,|\,\sigma\!,\,\omega\subset[m],\,\sigma{\cap}\omega=\emptyset\}$.
For $(\sigma\!,\,\omega)\in\XX_m$, define simplicial complex
$K_{\sigma\!,\omega}={\rm link}_K\sigma|_\omega=\{\tau\subset\omega\,|\,\sigma{\cup}\tau\in K\}$
(so $K_{\sigma\!,\omega}=\{\,\}$ if $\sigma\notin K$ or $K=\{\,\}$).
Then for any $(\sigma\!,\,\omega)\in\XX_m$ such that $\omega\neq\emptyset$,
$$(K_{\sigma\!,\,\omega})^\circ=(K^\circ)_{\tilde\sigma\!,\,\omega},\quad\w\sigma=[m]{\setminus}(\sigma{\cup}\omega),$$
where $(K_{\sigma\!,\,\omega})^\circ$ is the dual of
$K_{\sigma\!,\,\omega}$ relative to $\omega$.
}\vspace{3mm}

{\it Proof}\, Suppose $\sigma\in K$. Then
$$\begin{array}{l}
\quad (K^\circ)_{\tilde\sigma\!,\,\omega}\vspace{1mm}\\
=\{\eta\,\,|\,\,\eta\subset\omega,\,\,[m]{\setminus}(\w\sigma{\cup}\eta)=\sigma\!\cup\!(\omega\!\setminus\!\eta)\notin K\,\}\vspace{1mm}\\
=\{\omega\!\setminus\!\tau\,\,|\,\,\tau\subset\omega,\,\,\sigma\!\cup\!\tau\notin K\,\}
\quad(\tau=\omega\!\setminus\!\eta)\vspace{1mm}\\
=(K_{\sigma\!,\,\omega})^\circ.
\end{array}$$

If $\sigma\notin K$, then $(K^\circ)_{\tilde\sigma\!,\omega}=\Delta\!^\omega=(K_{\sigma\!,\omega})^\circ$.
\hfill $\Box$\vspace{3mm}

A sequence of simplicial pairs $(\underline{X},\underline{A})=\{(X_k,A_k)\}_{k=1}^m$ in this paper means that
the vertex set of $X_k$ is a subset of $[n_k]$ ($n_k>0$) which is the subset
$$\{s_{k-1}{+}1,s_{k-1}{+}2,\cdots,s_{k-1}{+}n_k\}\quad(s_k=n_1{+}{\cdots}{+}n_k,\,\,s_{0}=0)$$
of $[n]$ with $n=n_1{+}{\cdots}{+}n_m$.

For simplicial complexes $Y_1,{\cdots},Y_m$ such that the vertex set of $Y_k$ is a subset of $[n_k]$,
the union simplicial complex is
$$Y_1*{\cdots}*Y_m=\{\sigma\!\subset\![n]\,|\,\sigma{\cap}[n_k]\in Y_k\,\,{\rm for}\,\,k\!=\!1,{\cdots},m\}.$$

{\bf Definition 2.7} Let $K$ be a simplicial complex on $[m]$ and
$(\underline{X},\underline{A})$ be as above.
The {\it polyhedral join complex} ${\cal S}(K;\underline{X},\underline{A})$
is the simplicial complex on $[n]$ defined as follows. For a subset $\tau\subset[m]$,
define
$$S(\tau)= Y_1*\cdots *Y_m,\quad Y_k=\left\{\begin{array}{cl}
X_k&{\rm if}\,\,k\in \tau, \vspace{1mm}\\
A_k&{\rm if}\,\,k\not\in \tau.
\end{array}
\right.$$
Then ${\cal S}(K;\underline{X},\underline{A})=\cup_{\tau\in K}S(\tau)$.
Void complex $\{\,\}$ is allowed in a simplicial pair and $\{\,\}*X=\{\,\}$ for all $X$.
Define ${\cal Z}(\{\,\};\underline{X},\underline{A})=\{\,\}$.
\vspace{3mm}

The name polyhedral join comes from Definition 4.2 of \cite{AY}.
\vspace{3mm}

{\bf Example 2.8} For ${\cal S}(K;\underline{X},\underline{A})$,
let $S=\{k\,|\,A_k=\{\,\}\}$. Then
$${\cal S}(K;\underline{X},\underline{A})={\cal S}({\rm link}_KS;\underline{X'},\underline{A'}){*}(*_{k\in S}X_k),$$
where $(\underline{X'},\underline{A'})=\{(X_k,A_k)\}_{k\notin S}$.
\vspace{3mm}

{\bf Theorem 2.9} {\it Let ${\cal Z}(K;\underline{Y},\underline{B})$,
$(\underline{Y},\underline{B})=\{(Y_k,B_k)\}_{k=1}^m$ be the polyhedral product space
defined as follows. For each $k$, $(Y_k,B_k)$ is a pair of polyhedral product spaces given by
$$(Y_k,B_k)=\big({\cal Z}(X_k;\underline{U_k},\underline{C_k}),\,{\cal Z}(A_k;\underline{U_k},\underline{C_k})\big),\,\,
(\underline{U_k},\underline{C_k})=\{(U_i,C_i)\}_{i\,=\,s_{k-1}{+}1}^{s_k},$$
where $(X_k,A_k)$ is a simplicial pair on $[n_k]$.
Then
$${\cal Z}(K;\underline{Y},\underline{B})={\cal Z}({\cal S}(K;\underline{X},\underline{A});\underline{U},\underline{C}),$$
where $(\underline{U}{,}\underline{C})=\{(U_k{,}C_k)\}_{k=1}^n$,
$n=n_1{+}{\cdots}{+}n_m$.
\vspace{2mm}

Proof}\, If $K=\{\,\}$ or $U_k=\emptyset$ for some $k$ (this implies $Y_l=\emptyset$ for some $l$), then
${\cal Z}(K;\underline{Y},\underline{B})={\cal Z}({\cal S}(K;\underline{X},\underline{A});\underline{U},\underline{C})=\emptyset$.
So we suppose $K\neq\{\,\}$ and $U_k\neq\emptyset$ for all $k$ in the remaining part of proof.

We first prove the case $C_k\neq\emptyset$ for all $k$.
If $X_k=\{\,\}$ for some $k$, then $Y_k=\emptyset$ and ${\cal S}(K;\underline{X},\underline{A})=\{\,\}$.
So ${\cal Z}(K;\underline{Y},\underline{B})
={\cal Z}({\cal S}(K;\underline{X},\underline{A});\underline{U},\underline{C})=\emptyset$.
Suppose $X_k\neq\{\,\}$ for all $k$.
Let $S=\{k\,|\,A_k=\{\,\}\}=\{k\,|\,B_k=\emptyset\}$. If $S\not\in K$, then ${\rm link}_KS=\{\,\}$.
From Example~2.2 and Example 2.8 we have
${\cal Z}(K;\underline{Y},\underline{B})={\cal Z}({\cal S}(K;\underline{X},\underline{A});\underline{U},\underline{C})=\emptyset$.
Suppose $S\in K$.
Let $Z^\tau_k=Y_k$ if $k\in\tau$ and $Z^\tau_k=B_k$ if $k\notin\tau$.
For $\tau_k\subset[n_k]$, $W^{\tau_k}_t=U_k$ if $t\in\tau_k$ and $W^{\tau_k}_t=C_k$ if $k\in[n_k]{\setminus}\tau_k$. Then
$$\begin{array}{l}
\quad {\cal Z}(K;\underline{Y},\underline{B})\vspace{2mm}\\
=\cup_{\tau\in K} Z^\tau_1\times{\cdots}\times Z^\tau_m\vspace{2mm}\\
=\cup_{\tau,\tau_1,{\cdots},\tau_m}\,(W^{\tau_1}_1{\times}{\cdots}{\times}W^{\tau_1}_{n_1}){\times}{\cdots}{\times}
(W^{\tau_m}_{s_{m-1}+1}{\times}{\cdots}{\times}W^{\tau_m}_n)\vspace{2mm}\\
=\cup_{(\tau_1,{\cdots},\tau_m)\in{\cal S}(K;\underline{X},\underline{A})}\,(W^{\tau_1}_1{\times}{\cdots}{\times}W^{\tau_1}_{n_1}){\times}{\cdots}{\times}
(W^{\tau_m}_{s_{m-1}+1}{\times}{\cdots}{\times}W^{\tau_m}_n)\vspace{2mm}\\
={\cal Z}({\cal S}(K;\underline{X},\underline{A});\underline{U},\underline{C}),
\end{array}$$
where $\tau,\tau_1,{\cdots},\tau_m$ are taken over all subsets such that $\tau\in K$, $\tau_k\in X_k$ if $k\in\tau$
and $\tau_k\in A_k$ if $k\notin\tau$.

Now we prove the case $\sigma=\{k\,|\, C_k=\emptyset\}\neq\emptyset$.
Let $\sigma_k=\sigma{\cap}[n_k]$. Then from Example 2.2 we have
$$\begin{array}{l}
\,Y_k=Y'_k\times(\times_{k\in\sigma_k}U_k),\quad
\,Y'_k={\cal Z}({\rm link}_{X_k}\sigma_k;\underline{U'_k},\underline{C'_k}),\vspace{2mm}\\
B_k=B'_k\times(\times_{k\in\sigma_k}U_k),\quad
B'_k={\cal Z}({\rm link}_{A_k}\sigma_k;\underline{U'_k},\underline{C'_k}),\vspace{2mm}\\
{\cal Z}(K;\underline{Y},\underline{B})={\cal Z}(K;\underline{Y'},\underline{B'})\times(\times_{k\in\sigma}U_k),
\end{array}$$
where $(\underline{U'_k},\underline{C'_k})=\{(U_k,C_k)\}_{k\in[n_k]\setminus\sigma_k}$
and $(\underline{Y'},\underline{B'})=\{(Y'_k,B'_k)\}_{k\notin\sigma}$.
Then
$$\begin{array}{l}
\quad{\cal Z}(K;\underline{Y},\underline{B})\vspace{2mm}\\
={\cal Z}(K;\underline{Y'},\underline{B'})\times(\times_{k\in\sigma}U_k)\vspace{2mm}\\
={\cal Z}({\cal S}(K;{\rm link}_{(\underline{X},\underline{A})}\sigma);\underline{U'},\underline{C'})\times(\times_{k\in\sigma}U_k)\vspace{2mm}\\
={\cal Z}({\rm link}_{{\cal S}(K;\underline{X},\underline{A})}\sigma;\underline{U'},\underline{C'})\times(\times_{k\in\sigma}U_k)
\quad ({\rm by}\,{\rm Theorem\,2.10})\vspace{2mm}\\
={\cal Z}({\cal S}(K;
\underline{X},\underline{A});\underline{U},\underline{C}),
\end{array}$$
where $(\underline{U'},\underline{C'})=\{(U_k,C_k)\}_{k\notin\sigma}$.
\hfill$\Box$\vspace{3mm}

The above theorem is a trivial extension of Proposition 5.1 of \cite{AY}.
With this theorem we see that to compute the cohomology algebra of ${\cal Z}(K;\underline{Y},\underline{B})$,
we have to compute the universal algebra of ${\cal S}(K;\underline{X},\underline{A})$, which is the central work
of this paper.
\vspace{3mm}

{\bf Theorem 2.10} {\it \, For the ${\cal S}(K;\underline{X},\underline{A})$ in Definition 2.7 and
$(\sigma\!,\,\omega)\in\XX_n$ ( the simplicial complex $(-)_{\sigma\!,\,\omega}$ is as defined in Theorem 2.6),
$${\cal S}(K;\underline{X},\underline{A})_{\sigma,\,\omega}
={\cal S}(K;\underline{X}\,_{\sigma,\,\omega},\underline{A}\,_{\sigma,\,\omega}),$$
where $(\underline{X}\,_{\sigma,\,\omega},\underline{A}\,_{\sigma,\,\omega})\!=\!
\{(\,(X_k)_{\sigma_k,\,\omega_k},(A_k)_{\sigma_k,\,\omega_k})\}_{k=1}^m$,
$\sigma_k\!=\!\sigma{\cap}[n_k]$, $\omega_k\!=\!\omega{\cap}[n_k]$.

For $\sigma\subset[n]$, $\sigma_k=\sigma{\cap}[n_k]$ (the links as defined in conventions),
$${\rm link}_{{\cal S}(K;\underline{X},\underline{A})}\sigma=
{\cal S}(K;{\rm link}_{(\underline{X},\underline{A})}\sigma),$$
where ${\rm link}_{(\underline{X},\underline{A})}\sigma{=}\{({\rm link}_{X_k}\sigma_k,{\rm link}_{A_k}\sigma_k)\}_{k=1}^m$.
Precisely, if there is some $\sigma_k\notin X_k$,
then $\sigma\notin{\cal S}(K;\underline{X},\underline{A})$ and so
${\rm link}_{{\cal S}(K;\underline{X},\underline{A})}\sigma=\{\,\}$. Suppose
$\sigma_k\in X_k$ for all $k$. Let $\hat\sigma=\{k\,|\,\sigma_k\notin A_k\}$,
$\sigma'=\cup_{k\notin\hat\sigma}\,\sigma_k$, $n'=\Sigma_{k\in\hat\sigma}\,(n_k{-}|\sigma_k|)$. Then
$${\rm link}_{{\cal S}(K;\underline{X},\underline{A})}\sigma=
{\cal S}({\rm link}_K\hat\sigma;{\rm link}_{(\underline{X},\underline{A})'}\sigma')*\Delta\!^{[n']}
,$$
where ${\rm link}_{(\underline{X},\underline{A})'}\sigma'{=}\{({\rm link}_{X_k}\sigma_k,{\rm link}_{A_k}\sigma_k)\}_{k\notin\hat\sigma}$.
\vspace{2mm}

Proof}\, Let $Y^\tau_k=X_k$ if $k\in\tau$ and $Y^\tau_k=A_k$ if $k\notin\tau$ ($\{\,\}$ is allowed).
Then
\begin{eqnarray*}&&{\cal S}(K;\underline{X},\underline{A})_{\sigma\!,\,\omega}\vspace{1mm}\\
&=&\cup_{\tau\in K}(Y^\tau_1*{\cdots}*Y^\tau_m)_{\sigma\!,\,\omega}\vspace{1mm}\\
&=&\cup_{\tau\in K}(Y^\tau_1)_{\sigma_1,\,\omega_1}*\cdots*(Y^\tau_m)_{\sigma_m,\,\omega_m}\vspace{1mm}\\
&=&{\cal S}(K;\underline{X}\,_{\sigma,\,\omega},\underline{A}\,_{\sigma,\,\omega}).
\end{eqnarray*}

For a simplicial  complex $L$ on $[l]$ and $\sigma\subset[l]$, ${\rm link}_L\sigma=L_{\sigma,[l]\setminus\sigma}$.
So the second equality of the theorem is the special case of the first equality for $\omega=[n]{\setminus}\sigma$.
The third equality holds by Example 2.8.
\vspace{3mm}\hfill$\Box$

The dual of ${\cal S}(K;\underline{X},\underline{A})$ relative to $[n]$ is in general not a polyhedral join complex.
But the dual of the following type is.
\vspace{3mm}

{\bf Definition 2.11} The {\it composition complex} ${\cal S}(K;L_1,{\cdots},L_m)$ is the polyhedral join complex ${\cal S}(K;\underline{X},\underline{A})$
such that each $(X_k,A_k)=(\Delta\!^{[n_k]},L_k)$.
\vspace{3mm}

The name composition complex comes from Definition 4.5 of \cite{AY}.
\vspace{3mm}

{\bf Theorem 2.12} {\it Let ${\cal S}(K;L_1,{\cdots},L_m)^\circ$ be the dual of ${\cal S}(K;L_1,{\cdots},L_m)$ relative to $[n]$.
Then
$${\cal S}(K;L_1,{\cdots},L_m)^\circ={\cal S}(K^\circ;L_1^\circ,{\cdots},L_m^\circ),$$
where $K^\circ$ is the dual of $K$ relative to $[m]$ and $L_k^\circ$ is the dual  of $L_k$ relative to $[n_k]$.
So if $K$ and all $L_k$ are self dual ($X\cong X^\circ$ relative to its non-empty vertex set),
then ${\cal S}(K;L_1,{\cdots},L_m)$ is self dual.
\vspace{2mm}

Proof}\, For $\sigma\subset[m]$ but $\sigma\neq[m]$ ($\sigma=\emptyset$, $L_k=\Delta\!^{[n_k]}$ or $\{\,\}$ are allowed),
$$\begin{array}{l}
\quad {\cal S}(\Delta\!^\sigma;L_1,{\cdots},L_m)^\circ\vspace{2mm}\\
=\{[n]\!\setminus\!\tau\,|\,\tau\in\cup_{j\notin\sigma}\Delta\!^{n_1}*{\cdots}*(\Delta\!^{n_j}{\setminus}L_j)*{\cdots}*\Delta\!^{n_m}\}\vspace{2mm}\\
=\cup_{j\notin\sigma}\,\Delta\!^{n_1}*{\cdots}*L_j^\circ*{\cdots}*\Delta\!^{n_m}\vspace{2mm}\\
={\cal S}((\Delta\!^\sigma)^\circ;L_1^\circ,{\cdots},L_m^\circ),
\end{array}$$
So for $K\neq[m]$ or $\{\,\}$,
$$\begin{array}{l}
\quad {\cal S}(K;L_1,{\cdots},L_m)^\circ\vspace{2mm}\\
=(\cup_{\sigma\in K}{\cal S}((\Delta\!^\sigma);L_1,{\cdots},L_m))^\circ\vspace{2mm}\\
={\cal S}(\cap_{\sigma\in K}(\Delta\!^\sigma)^\circ;L_1^\circ,{\cdots},L_m^\circ)\vspace{2mm}\\
={\cal S}(K^\circ;L_1^\circ,{\cdots},L_m^\circ).
\end{array}$$

For $K=\Delta\!^{[m]}$ or $\{\,\}$, the equality holds naturally.
\hfill$\Box$\vspace{3mm}

\section{Homology and Cohomology Group}

This is a paper following \cite{Z}.
All the basic definitions such as indexed groups and (co)chain complexes, diagonal tensor product,
etc., are as in \cite{Z}.

In this section, we prove that the (co)homology group of the polyhedral product inclusion complex
$C_*(K;\underline{\vartheta})$ is the diagonal tensor product of the total (co)homology group of $K$
and the character group of the induced homomorphism $\underline{\theta}$ ($\underline{\theta}^\circ$) in Theorem 3.7.
\vspace{3mm}

{\bf Conventions} In this paper, a group $A_*^\Lambda=\oplus_{\alpha\in\Lambda}\,A_*^\alpha$ indexed by $\Lambda$
is simply denoted by $A_*$ when there is no confusion.
So is the (co)chain complex case.
The diagonal tensor product $A_*^\Lambda{\otimes}_\Lambda B_*^\Lambda$ in \cite{Z} is simply denoted by
$A_*^\Lambda\,\widehat\otimes\,B_*^\Lambda$ in this paper (the index set $\Lambda$ can not be omitted in this case).
\vspace{3mm}

{\bf Definition 3.1} Let $A_*=\oplus_{\alpha\in\Lambda}A_*^\alpha$, $B_*=\oplus_{\alpha\in\Lambda}B_*^\alpha$
be two groups indexed by the same set $\Lambda$.
An {\it indexed group homomorphism} $f\colon A_*\to B_*$ is the direct sum
$f=\oplus_{\alpha\in\Lambda}f_\alpha$ such that each $f_\alpha\colon A_*^\alpha\to B_*^{\alpha}$
is a graded group homomorphism. Define groups indexed by $\Lambda$ as follows.
$$\begin{array}{ll}
{\rm ker}\,f=\oplus_{\alpha\in\Lambda}{\rm ker}\,f_\alpha,&{\rm coker}\,f=\oplus_{\alpha\in\Lambda}{\rm coker}\,f_\alpha,\vspace{1mm}\\
{\rm im}\,f\,=\oplus_{\alpha\in\Lambda}\,{\rm im}\,f_\alpha,&{\rm coim}\,f\,=\oplus_{\alpha\in\Lambda}\,{\rm coim}\,f_\alpha.
\end{array}$$

For indexed group homomorphism $f=\oplus_{\alpha\in\Lambda}\,f_\alpha$ and $g=\oplus_{\beta\in\Gamma}\,g_\beta$,
their tensor product $f{\otimes}g$ is naturally an indexed group homomorphism with
$f{\otimes}g=\oplus_{(\alpha,\beta)\in\Lambda{\times}\Gamma}f_\alpha{\otimes}g_\beta$.

For indexed group homomorphism $f=\oplus_{\alpha\in\Lambda}\,f_\alpha$ and $g=\oplus_{\alpha\in\Lambda}\,g_\alpha$ indexed by
the same set,
their diagonal tensor product $f\widehat\otimes g$ is the indexed group homomorphism
$f\widehat\otimes g=\oplus_{\alpha\in\Lambda}f_\alpha{\otimes}g_\alpha$.

Similarly, we have the definition of {\it indexed (co)chain homomorphism}
by replacing the indexed groups in the above definition by indexed (co)chain complexes.
\vspace{3mm}

{\bf Definition 3.2} An indexed group homomorphism $\theta\colon U_*\to V_*$
is called a {\it split homomorphism} if ${\rm ker}\,\theta$, ${\rm coker}\,\theta$ and ${\rm im}\,\theta$ are all free groups.

An indexed chain homomorphism
$\vartheta\colon (C_*,d)\to(D_*,d)$ with induced homology group homomorphism
$\theta\colon U_*\!=\!H_*(C_*)\to V_*\!=\!H_*(D_*)$
 is called a {\it split inclusion} if $C_*$ is a chain subcomplex of the free complex $D_*$ and
$\theta$ is a split homomorphism.
\vspace{3mm}

{\bf Definition 3.3} Let $\theta\colon U_*\to V_*$ be a split homomorphism
with dual homomorphism $\theta^\circ\colon V^*\to U^*$.
The index set $\XX=\{\,{\scriptstyle(\emptyset,\emptyset)},{\scriptstyle(\emptyset,\{1\})},
{\scriptstyle(\{1\},\emptyset)}\,\}$ and
$\RR=\{\,{\scriptstyle(\emptyset,\emptyset)},{\scriptstyle(\emptyset,\{1\})}\,\}\subset\XX$.
$\SS$ may be taken to be either $\XX$ or $\RR$.

The indexed groups $H_*^{\SS}(\theta)=\oplus_{s\in\SS}\,H_*^{s}(\theta)$
and its dual groups $H^{\,*}_{\SS}(\theta^\circ)=\oplus_{s\in\SS}\,H^*_{s}(\theta^\circ)$
are given by
$$H_*^s(\theta)=\left\{\begin{array}{ll}
{\rm coker}\,\theta&{\rm if}\,\,s={\scriptstyle(\{1\},\emptyset)},\\
{\rm ker}\,\theta&{\rm if}\,\,s={\scriptstyle(\emptyset,\{1\})},\\
{\rm im}\,\theta&{\rm if}\,\,s={\scriptstyle(\emptyset,\emptyset)},
\end{array}\right.$$
$$H^*_s(\theta^\circ)=\left\{\begin{array}{ll}
{\rm ker}\,\theta^\circ&{\rm if}\,\,s={\scriptstyle(\{1\},\emptyset)},\\
{\rm coker}\,\theta^\circ&{\rm if}\,\,s={\scriptstyle(\emptyset,\{1\})},\\
{\rm im}\,\theta^\circ&{\rm if}\,\,s={\scriptstyle(\emptyset,\emptyset)},
\end{array}\right.
$$
$H_*^{\SS}(\theta)$ ($H^{\,*}_{\SS}(\theta^\circ)$) is called the {\it (right {\rm for} $\RR$) character group} of $\theta$ ($\theta^\circ$).

The indexed chain complexes $(C_*^{\SS}(\theta),d)=\oplus_{s\in\SS}\,(C_*^{s}(\theta),d)$
and its dual cochain complexes $(C^{\,*}_{\!\SS}(\theta^\circ),\delta)=\oplus_{s\in\SS}\,(C^*_{s}(\theta^\circ),\delta)$
are given by
$$C_*^s(\theta)=\left\{\begin{array}{ll}
{\rm coker}\,\theta&{\rm if}\,\,s={\scriptstyle(\{1\},\emptyset)},\\
{\rm ker}\,\theta{\oplus}\Sigma{\rm ker}\,\theta&{\rm if}\,\,s={\scriptstyle(\emptyset,\{1\})},\\
{\rm im}\,\theta&{\rm if}\,\,s={\scriptstyle(\emptyset,\emptyset)},
\end{array}\right.\quad\quad$$
$$C^*_s(\theta^\circ)=\left\{\begin{array}{ll}
{\rm ker}\,\theta^\circ&{\rm if}\,\,s={\scriptstyle(\{1\},\emptyset)},\\
{\rm coker}\,\theta^\circ{\oplus}\Sigma{\rm coker}\,\theta^\circ&{\rm if}\,\,s={\scriptstyle(\emptyset,\{1\})},\\
{\rm im}\,\theta^\circ&{\rm if}\,\,s={\scriptstyle(\emptyset,\emptyset)},
\end{array}\right.$$
where $d$ is trivial on $C_*^{\emptyset,\emptyset}(\theta)$ and $C_*^{\{1\},\emptyset}(\theta)$
and is the desuspension isomorphism on $C_*^{\emptyset,\{1\}}(\theta)$.
$C_*^{\SS}(\theta)$ ($C^{\,*}_{\SS}(\theta^\circ)$) is called the {\it (right {\rm for} $\RR$) character complex} of $\theta$ ($\theta^\circ$).

Let $\theta=\oplus_{\alpha\in\Lambda}\theta_{\alpha}$.
Then $H_*^{\SS}\!(\theta)$ is also a group indexed by $\Lambda$
and so denoted by $H_*^{\SS}\!(\theta)=H_*^{\SS\!;\Lambda}(\theta)=\oplus_{s\in\SS\!,\alpha\in\Lambda}\,
H_*^{s;\alpha}(\theta)$ with
$$H_*^{s;\alpha}(\theta)=\left\{\begin{array}{ll}
{\rm coker}\,\theta_\alpha&{\rm if}\,\,s={\scriptstyle(\{1\},\emptyset)},\\
{\rm ker}\,\theta_\alpha&{\rm if}\,\,s={\scriptstyle(\emptyset,\{1\})},\\
{\rm im}\,\theta_\alpha&{\rm if}\,\,s={\scriptstyle(\emptyset,\emptyset)}.
\end{array}\right.$$
Other cases are similar.
\vspace{3mm}

{\bf Theorem 3.4}\, {\it For a split inclusion $\vartheta\colon(C_*,d)\to(D_*,d)$
with induced homology homomorphism $\theta\colon U_*\to V_*$,
there are quotient chain homotopy equivalences  $q$ and $q'$
satisfying the following commutative diagram
$$\begin{array}{ccc}
(C_*,d)&\stackrel{q'}{\longrightarrow}&(U_*,d)\,\vspace{1mm}\\
^\vartheta\downarrow\,\,&&^{\vartheta'}\downarrow\quad\\
(D_*,d)&\stackrel{q}{\longrightarrow}&(C_*^{\XX}(\theta),d),
\end{array}
$$
where $\vartheta'$ is the inclusion by identifying $U_*={\rm ker}\,\theta\oplus{\rm coim}\,\theta$ with
${\rm ker}\,\theta\oplus{\rm im}\,\theta\subset C_*^{\XX}(\theta)$ (as in Definition 3.3).

There are also isomorphisms $\phi$ and $\phi'$ of chain complexes indexed by $\XX$
satisfying the following commutative diagram
$$\begin{array}{ccc}
(U_*,d)&\stackrel{\phi'}{\longrightarrow}&S_*^\XX\,\widehat\otimes\, H_*^{\XX}\!(\theta)\vspace{1mm}\\
^{\vartheta'}\downarrow\,\,&&^{i\,\widehat\otimes\, 1}\downarrow\quad\quad\quad\\
(C_*^{\XX}\!(\theta),d)&\stackrel{\phi}{\longrightarrow}&(T_*^\XX\,\widehat\otimes\, H_*^{\XX}\!(\theta),d)
\end{array}
$$
where $S_*^\XX$ and $T_*^\XX$ are as defined in} \cite{Z}, {\it $1$ is the identity and $i$ is the inclusion.
Precisely, let $\Bbb Z(x_1,{\cdots},x_n)$ be the free abelian group generated by $x_1,{\cdots},x_n$, then
$(T_*^{\emptyset,\emptyset},d)=\zz(\eta),\,\,(T_*^{\emptyset,\{1\}},d)=(\zz(\beta,\gamma),d),
(T_*^{\{1\},\emptyset},d)=\zz(\alpha)$,
$|\alpha|=|\beta|{-}1=|\gamma|=|\eta|=0$, $d\beta=\gamma$,
$S_*^\XX=\Bbb Z(\gamma,\eta)$.

If $\theta$ is an epimorphism, then
$H_*^{\XX}\!(\theta)=H_*^{\RR}(\theta)=U_*$ by identifying ${\rm im}\,\theta$
with ${\rm coim}\,\theta$ and so  all $\XX$ is replaced by $\RR$.
\vspace{2mm}

Proof}\, Take a representative $a_i$ in $C_*$ for every generator of ${\rm ker}\,\theta$
and let $\overline a _i\in D_*$ be any element such that $d\overline a _i= a_i$.
Take a representative $b_j$ in $C_*$ for every generator of ${\rm im}\,\theta$.
Take a representative $c_k$ in $D_*$ for every generator of ${\rm coker}\,\theta$.
So we may regard $U_*$ as the chain subcomplex of $C_*$ freely generated by all $a_i$'s and $b_j$'s
and regard $(C_*^{\XX\!}(\theta),d)$ as the chain subcomplex of $D_*$ freely generated by all
$a_i$'s, $\overline a _i$'s, $b_j$'s and $c_k$'s.
Then we have the following commutative diagram of short exact sequences of chain complexes
$$\begin{array}{ccccccc}
 0\to& U_*&\stackrel{i}{\longrightarrow}& C_*&\stackrel{j}{\longrightarrow}&C_*/U_*&\to0\vspace{1mm}\\
 &^{\vartheta'}\downarrow\quad&&^\vartheta\downarrow\,\quad&&\downarrow&\\
0\to& C_*^{\XX\!}(\theta)&\stackrel{i}{\longrightarrow}& D_*&\stackrel{j}{\longrightarrow}&D_*/C_*^{\XX\!}(\theta)&\to0.
  \end{array}
$$
Since all the complexes are free, the two $i$'s have group homomorphism inverse. $H_*(C_*/U_*)=0$ and $H_*(D_*/C_*^{\XX\!}(\theta))=0$
imply that the inverse of $i$'s are complex homomorphisms. So we may take $q,q'$ to be the inverse of $i$'s.

$\phi$ is defined as shown in the following table.
\begin{center}
\begin{tabular}{|c|c|c|c|c|}
\hline
{\rule[-2mm]{0mm}{6mm}$\quad x\in$}&${\rm coker}\,\theta$&$\Sigma\,{\rm ker}\,\theta$
&${\rm ker}\,\theta$&${\rm im}\,\theta$\\
\hline
{\rule[-2mm]{0mm}{7mm}$\phi(x)=$}&$\alpha\,\widehat\otimes\,x$&$\beta\,\widehat\otimes\,dx$
&$\gamma\,\widehat\otimes\,x$&$\eta\,\widehat\otimes\,x$\\
\hline
\end{tabular}
\end{center}
\hfill$\Box$\vspace{3mm}

Compare the proof of the above theorem with that of Theorem 2.3 of \cite{Z}.
The work of this and the next section is just to generalize all definitions for
graded groups $A_*$ ($A_*^\XX$) in \cite{Z} to definitions for indexed groups $A_*=\oplus_{\alpha\in\Lambda}\,A_*^\alpha$
($A_*^\XX=\oplus_{\alpha\in\Lambda}A_*^{\XX;\alpha}$). Then all the proofs in this paper are naturally obtained from \cite{Z} by
replacing all the graded groups ((co)chain complexes) by indexed groups ((co)chain complexes). So we omit
the proof when the analogue in \cite{Z} is given.
\vspace{3mm}

{\bf Definition 3.5} For $k=1,{\cdots},m$,
let $\vartheta_k\colon\! ((C_k)_*,d)\to((D_k)_*,d)$ be a split inclusion with induced
homology group homomorphism $\theta_k\colon\! (U_k)_*\to (V_k)_*$.
Denote $\underline{\vartheta}=\{\vartheta_k\}_{k=1}^m$, $\underline{\theta}=\{\theta_k\}_{k=1}^m$
and their dual $\underline{\vartheta}^\circ=\{\vartheta_k^\circ\}_{k=1}^m$, $\underline{\theta}^\circ=\{\theta_k^\circ\}_{k=1}^m$.
The index set $\SS=\XX$ or $\RR$.

The indexed group $H_*^{\SS_m}(\underline{\theta})$ and its dual group $H^{\,*}_{\!\SS_m}(\underline{\theta}^\circ)$
are given by
$$H_*^{\SS_m}(\underline{\theta})=
H_*^{\SS}(\theta_1)\otimes{\cdots}\otimes H_*^{\SS}(\theta_m),\,\,
H^{\,*}_{\!\SS_m}(\underline{\theta}^\circ)=
H^{\,*}_{\!\SS}(\theta_1^\circ)\otimes{\cdots}\otimes H^{\,*}_{\!\SS}(\theta_m^\circ).$$
Since $\SS_m=\SS{\times}{\cdots}{\times}\SS$ ($m$-fold) by the 1-1 correspondence
$$(\sigma\!,\,\omega)\to(s_1,{\cdots},s_m),\,\,\sigma=\{k\,|\,s_k={\scriptstyle(\{1\},\emptyset)}\},\,\,
\omega=\{k\,|\,s_k={\scriptstyle(\emptyset,\{1\})}\}.$$
we may write $H_*^{\SS_m}(\underline{\theta})=\oplus_{(\sigma\!,\,\omega)\in\SS_m}H_*^{\sigma\!,\,\omega}(\underline{\theta})$,
$H^{\,*}_{\!\SS_m}(\underline{\theta}^\circ)=\oplus_{(\sigma\!,\,\omega)\in\SS_m}H^*_{\sigma\!,\,\omega}(\underline{\theta}^\circ)$.
Then by definition,
$$H_*^{\sigma\!,\,\omega}(\underline{\theta})
= H_1{\otimes}{\cdots}{\otimes}H_m,\quad
H_k=\left\{\begin{array}{cl}
{\rm coker}\,\theta_{k}&{\rm if}\, k\in\sigma,\vspace{1mm}\\
{\rm ker}\,\theta_{k}&{\rm if}\, k\in\omega,\vspace{1mm}\\
{\rm im}\,\theta_{k}&{\rm otherwise},
\end{array}\right.$$
$$H^*_{\sigma\!,\,\omega}(\underline{\theta}^\circ)
= H^1{\otimes}{\cdots}{\otimes}H^m,\quad
H^k=\left\{\begin{array}{cl}
{\rm ker}\,\theta_{k}^\circ&{\rm if}\, k\in\sigma,\vspace{1mm}\\
{\rm coker}\theta_{k}^\circ&{\rm if}\, k\in\omega,\vspace{1mm}\\
{\rm im}\,\theta_{k}^\circ&{\rm otherwise}.
\end{array}\right.$$
$H_*^{\SS_m}(\underline{\theta})$ ($H^*_{\SS_m}(\underline{\theta}^\circ)$) is called the {\it (right {\rm for} $\RR$) character group} of
$\underline{\theta}$ ($\underline{\theta}^\circ$).

The indexed chain complex $(C_*^{\SS_m}(\underline{\theta}),d)$ and its dual cochain complex
$(C^{\,*}_{\!\SS_m}(\underline{\theta}^\circ),\delta)$ are given by
$$C_*^{\SS_m}(\underline{\theta})=C_*^{\SS}\!(\theta_1)\otimes{\cdots}\otimes C_*^{\SS}\!(\theta_m),\,\,
C^{\,*}_{\!\SS_m}(\underline{\theta}^\circ)=C^{\,*}_{\!\SS}(\theta_1^\circ)\otimes{\cdots}\otimes C^{\,*}_{\!\SS}(\theta_m^\circ).$$
$C_*^{\SS_m}(\underline{\theta})$ ($C^*_{\SS_m}(\underline{\theta}^\circ)$) is called the {\it (right {\rm for} $\RR$) character complex} of
$\underline{\theta}$ ($\underline{\theta}^\circ$).
\vspace{3mm}

{\bf Definition 3.6} Let $K$ be a simplicial complex on $[m]$
and everything else be as in Definition 3.5.
For $\SS=\XX$ or $\RR$, $T_*^\SS,S_*^\SS$ are as in Theorem 3.4.

The {\it total chain complex $(T_*^{\XX_m}(K),d)$} of $K$ is the chain subcomplex of $(T_*^{\XX_m}=T_*^\XX{\otimes}{\cdots}{\otimes}T_*^\XX,d)$
defined as follows.
For a subset $\tau$ of $[m]$, define
$$(U_*(\tau),d)=(U_1{\otimes}\cdots{\otimes}U_m,d),\quad (U_k,d)=\left\{\begin{array}{cl}
(T_*^\XX,d)&{\rm if}\,\,k\in \tau,\vspace{1mm}\\
S_*^\XX&{\rm if}\,\,k\not\in \tau,
\end{array}
\right.$$
Then $(T_*^{\XX_m}(K),d)=(+_{\tau\in K}\,U_*(\tau),d)$.
Define $(T_*^{\XX_m}(\{\}),d)=0$. The group $H_*^{\XX_m}(K)=H_*(T_*^{\XX_m}(K))$ is called {\it the
total homology group} of $K$.

So the  dual cochain complex $(T^*_{\XX_m}(K),\delta)$ of $(T_*^{\XX_m}(K),d)$
is a quotient complex of $(T^*_{\XX_m},\delta)$ and is called the {\it total cochain complex} of $K$.
The group $H^*_{\XX_m}(K)=H^*(T^*_{\XX_m}(K))$ is called {\it the
total cohomology group of $K$}.

The {\it polyhedral product character complex} $(C_*^{\XX_m}(K;\underline{\theta}),d)$ is the subcomplex of
$(C_*^{\XX_m}(\underline{\theta}),d)$ defined as follows.
For a subset $\tau$ of $[m]$, define
$$(H_*(\tau),d)=(H_1{\otimes}\cdots{\otimes}H_m,d),\quad (H_k,d)=\left\{\begin{array}{cl}
(C_*^{\XX}\!(\theta_k),d)&{\rm if}\,\,k\in \tau,\vspace{1mm}\\
((U_k)_*,d)&{\rm if}\,\,k\not\in \tau.
\end{array}\right.$$
Then $(C_*^{\XX_m}(K;\underline{\theta}),d)=(+_{\tau\in K}\,H_*(\tau),d)$. Define
$(C_*^{\XX_m}(\{\,\};\underline{\theta}),d)=0$.

So the dual cochain complex $(C^{\,*}_{\!\XX_m}(K;\underline{\theta}^\circ),\delta)$ of
$(C_*^{\XX_m}(K;\underline{\theta}),d)$ is a quotient complex of $(C^{\,*}_{\!\XX_m}(\underline{\theta}^\circ),\delta)$.

Replace $\XX$ by $\SS$ in the above definitions, we get the right analogues
({\it right total complex, right polyhedral product character complex,} etc.).

The {\it polyhedral product inclusion complex} $(C_*(K;\underline{\vartheta}),d)$ is the subcomplex of $((D_1)_*{\otimes}{\cdots}{\otimes}(D_m)_*,d)$
defined as follows.
For a subset $\tau$ of $[m]$, define
$$(E_*(\tau),d)=(E_1{\otimes}\cdots{\otimes}E_m,d),\quad (E_k,d)=\left\{\begin{array}{cl}
((D_k)_*,d)&{\rm if}\,\,k\in \tau,\vspace{1mm}\\
((C_k)_*,d)&{\rm if}\,\,k\not\in \tau.
\end{array}
\right.$$
Then $(C_*(K;\underline{\vartheta}),d)=(+_{\tau\in K}\,E_*(\tau),d)$.
Define $(C_*(\{\,\};\underline{\vartheta}),d)=0$.

So the  dual cochain complex $(C^*(K;\underline{\vartheta}^\circ),\delta)$ of $(C_*(K;\underline{\vartheta}),d)$
is a quotient complex of $((D_1)^*{\otimes}{\cdots}{\otimes}(D_m)^*,\delta)$.
\vspace{3mm}

By Theorem 4.7 of \cite{Z}, $H_*^{\sigma\!,\,\omega}(K)\cong\w H_{*-1}(K_{\sigma\!,\,\omega})$,
$H^*_{\sigma\!,\,\omega}(K)\cong\w H^{*-1}(K_{\sigma\!,\,\omega})$, where
$H_*^{\XX_m}(K)=\oplus_{(\sigma\!,\,\omega)\in\XX_m}\,H_*^{\sigma\!,\,\omega}(K)$,
$H^*_{\XX_m}(K)=\oplus_{(\sigma\!,\,\omega)\in\XX_m}\,H^*_{\sigma\!,\,\omega}(K)$.
\vspace{3mm}

{\bf Theorem 3.7}\, {\it For the $K$, $\underline{\vartheta}$ and $\underline{\theta}$ in Definition 3.5 and Definition 3.6,
there is a quotient chain homotopy equivalence (\,$\XX_m$ neglected)
$$\varphi_{(K;\underline{\vartheta})}\colon(C_*(K;\underline{\vartheta}),d)
\stackrel{\simeq}{\longrightarrow}
(C_*^{\XX_m}(K;\underline{\theta}),d)$$
and an isomorphism of chain complexes indexed by $\XX_m$
$$\phi_{(K;\underline{\theta})}\colon(C_*^{\XX_m}(K;\underline{\theta}),d)
\stackrel{\cong}{\longrightarrow}
(T_*^{\XX_m}(K)\,\widehat\otimes\,H_*^{\XX_m}(\underline{\theta}),d).$$
So we have (co)homology group isomorphisms
$$H_*(C_*(K;\underline{\vartheta}))\cong
H_*^{\XX_m}(K)\,\widehat\otimes\,H_*^{\XX_m}(\underline{\theta})
=H_*^{\XX_m}(K)\,\widehat\otimes\,\big(H_*^{\XX}\!(\theta_1){\otimes}{\cdots}{\otimes}H_*^{\XX}\!(\theta_m)\big),$$
$$H^*(C^*(K;\underline{\vartheta}^\circ))\cong
H^{\,*}_{\!\XX_m}(K)\,\widehat\otimes\,H^{\,*}_{\!\XX_m}(\underline{\theta}^\circ)
=H^*_{\!\XX_m}(K)\,\widehat\otimes\,\big(H^*_{\!\XX}(\theta_1^\circ){\otimes}{\cdots}{\otimes}H^*_{\!\XX}(\theta_m^\circ)\big).$$

If each $\theta_k$ is an epimorphism, then
$H_*^{\XX}\!(\theta_k)=H_*^{\RR}(\theta_k)=(U_k)_*$, $H^{\,*}_{\XX}(\theta_k^\circ)=H^{\,*}_{\RR}(\theta_k^\circ)=(U_k)^*$
and so all $\XX$ is replaced by $\RR$.
\vspace{2mm}

Proof}\, Similar to that of Theorem 2.6 and Theorem 4.6 of \cite{Z}.
\hfill$\Box$\vspace{3mm}

{\bf Definition 3.8} A polyhedral join complex ${\cal S}(K;\underline{X},\underline{A})$ is {\it homology split}
if the reduced simplicial homology homomorphism $$\iota_k\colon\w H_*(A_k)\to\w H_*(X_k)$$ induced by inclusion is split
for $k=1,{\cdots},m$.

A polyhedral join complex ${\cal S}(K;\underline{X},\underline{A})$ is {\it total homology split}
if the reduced simplicial homology homomorphism
$$\iota_{\sigma_k,\omega_k}\colon \w H_{*}((A_k)_{\sigma_k,\omega_k})
\to \w H_{*}((X_k)_{\sigma_k,\omega_k})\quad((-)_{\sigma\!,\,\omega}\,\,{\rm as\,\,in\,\,Theorem\, 2.6})$$ induced by inclusion is split
for all $(\sigma_k,\omega_k)\in\XX_{n_k}$.
\vspace{3mm}

{\bf Theorem 3.9}{\it\, For homology split ${\cal S}(K;\underline{X},\underline{A})$,
$$\begin{array}{l}
\quad\w H_{*-1}({\cal S}(K;\underline{X},\underline{A}))\vspace{2mm}\\
\cong H_*^{\XX_m}(K)\,\widehat\otimes\,H_*^{\XX_m}(\underline{X},\underline{A})\vspace{2mm}\\
=H_*^{\XX_m}(K)\,\widehat\otimes\,\big(H_*^{\XX}\!(X_1,A_1){\otimes}{\cdots}{\otimes}H_*^{\XX}\!(X_m,A_m)\big),
\end{array}$$
$$\begin{array}{l}
\quad\w H^{*-1}({\cal S}(K;\underline{X},\underline{A}))\vspace{2mm}\\
\cong H^*_{\XX_m}(K)\,\widehat\otimes\,H^*_{\XX_m}(\underline{X},\underline{A})\vspace{2mm}\\
=H^*_{\XX_m}(K)\,\widehat\otimes\,\big(H^{\,*}_{\!\XX}(X_1,A_1){\otimes}{\cdots}{\otimes}H^{\,*}_{\!\XX}(X_m,A_m)\big),
\end{array}$$
where $H_*^{\XX_m}(-)=\oplus_{(\sigma\!,\,\omega)\in\XX_m}\,H_*^{\sigma\!,\,\omega}(-)$,
$H^*_{\XX_m}(-)=\oplus_{(\sigma\!,\,\omega)\in\XX_m}\,H^*_{\sigma\!,\,\omega}(-)$ with
$$H_*^s(X_k,A_k)=\left\{\begin{array}{rl}
\Sigma{\rm coker}\,\iota_k&{\rm if}\, s={\scriptstyle (\{1\},\emptyset)},\vspace{1mm}\\
\Sigma{\rm ker}\,\iota_k&{\rm if}\,s={\scriptstyle (\emptyset,\{1\})},\vspace{1mm}\\
\Sigma{\rm im}\,\iota_k&{\rm if}\,s={\scriptstyle (\emptyset,\emptyset)},
\end{array}\right.$$
$$H^*_{s}(X_k,A_k)=\left\{\begin{array}{rl}
\Sigma{\rm ker}\,\iota_k^\circ&{\rm if}\, s={\scriptstyle (\{1\},\emptyset)},\vspace{1mm}\\
\Sigma{\rm coker}\,\iota_k^\circ&{\rm if}\,s={\scriptstyle (\emptyset,\{1\})},\vspace{1mm}\\
\Sigma{\rm im}\,\iota_k^\circ&{\rm if}\,s={\scriptstyle (\emptyset,\emptyset)},
\end{array}\right.$$
$$H_*^{\sigma\!,\,\omega}(\underline{X},\underline{A})
=H_1{\otimes}{\cdots}{\otimes}H_m,\quad
H_k=\left\{\begin{array}{rl}
\Sigma{\rm coker}\,\iota_k&{\rm if}\, k\in\sigma,\vspace{1mm}\\
\Sigma{\rm ker}\,\iota_k&{\rm if}\, k\in\omega,\vspace{1mm}\\
\Sigma{\rm im}\,\iota_k&{\rm otherwise},
\end{array}\right.$$
$$H^*_{\sigma\!,\,\omega}(\underline{X},\underline{A})
=H^1{\otimes}{\cdots}{\otimes}H^m,\quad
H^k=\left\{\begin{array}{rl}
\Sigma{\rm ker}\,\iota_k^\circ&{\rm if}\, k\in\sigma,\vspace{1mm}\\
\Sigma{\rm coker}\,\iota_k^\circ&{\rm if}\, k\in\omega,\vspace{1mm}\\
\Sigma{\rm im}\,\iota_k^\circ&{\rm otherwise},
\end{array}\right.$$
where $\iota_k$ is as in Definition 3.8 with dual $\iota_k^\circ$ and
$\Sigma$ means suspension.

If each $\iota_k$ is an epimorphism, then all $\XX$ is replaced by $\RR$ and we have
$H_*^\RR(X_k,A_k)=H_*(A_k)$, $H^*_\RR(X_k,A_k)=H^*(A_k)$.

If the reduced simplicial (co)homology is taken over a field, then the conclusion holds for all
polyhedral join complexes.
\vspace{2mm}

Proof}\, A corollary of Theorem 3.7 by
taking $\underline{\vartheta}=\{\vartheta_k\}_{k=1}^m$ with split inclusion
$\vartheta_k\colon(\Sigma\w C_*(A_k),d)\to(\Sigma\w C_*(X_k),d)$
the suspension reduced simplicial complex inclusion. Regard this graded group inclusion as
an indexed chain homomorphism such that the index set has only  one element.
Then $(C_*(K;\underline{\vartheta}),d)=(\Sigma\w C_*({\cal S}(K;\underline{X},\underline{A})),d)$
and $H_*^{\XX}\!(\theta_k)=H_*^{\XX}\!(X_k,A_k)$.
\hfill$\Box$\vspace{3mm}

{\bf Example 3.10} For ${\cal S}(K;\underline{X},\underline{A}))={\cal S}(K;L_1,{\cdots},L_m)$ such that all
$H_*(L_k)$ is free,
each $\iota_k\colon\w H_*(L_k)\to\w H_*(\Delta\!^{[n_k]})\,(=\!0)$ is an epimorphism.
By definition, $H_*^{\emptyset,\emptyset}(\Delta\!^{[n_k]},L_k)=0$,
$H_*^{\emptyset,\{1\}}(\Delta\!^{[n_k]},L_k)=\w H_{*-1}(L_k)$,
$H_*^{\emptyset,\omega}(\underline{X},\underline{A})=0$ if $\omega\neq[m]$. So
$$H_*^{\RR_m}(\underline{X},\underline{A})=H_*^{\emptyset,[m]}(\underline{X},\underline{A})=\w H_{*-1}(L_1){\otimes}{\cdots}{\otimes}\w H_{*-1}(L_m),$$
$$H_*^{\RR_m}(K)\,\widehat\otimes\,H_*^{\RR_m}(\underline{X},\underline{A})=\w H_{*-1}(K){\otimes}\w H_{*-1}(L_1){\cdots}{\otimes}\w H_{*-1}(L_m).$$
So by Theorem 3.9,
$$\w H_{*-1}({\cal S}(K;L_1,{\cdots},L_m))\cong\w H_{*-1}(K){\otimes}\w H_{*-1}(L_1){\cdots}{\otimes}\w H_{*-1}(L_m),$$
$$\w H^{*-1}({\cal S}(K;L_1,{\cdots},L_m))\cong\w H^{*-1}(K){\otimes}\w H^{*-1}(L_1){\cdots}{\otimes}\w H^{*-1}(L_m).$$

The above homology equality implies that
${\cal S}(K;L_1,{\cdots},L_m)$ is homology spherical ($\w H_{*}(-)\cong\Bbb Z$, so $\{\emptyset\}$
is homology spherical but $\{\,\}$ is not) if and only if $K$ and all $L_k$ are homology spherical.
By Theorem 2.10, for $\sigma\in\w K={\cal S}(K;L_1,{\cdots},L_m)$, $\sigma_k=\sigma{\cap}[n_k]$, $\hat\sigma=\{k\,|\,\sigma_k\notin L_k\}$,
$n'=\Sigma_{k\in\hat\sigma}(n_k{-}|\sigma_k|)$,
$${\rm link}_{\w K}\sigma={\cal S}({\rm link}_K\hat\sigma;{\rm link}_{L_{i_1}}\!\sigma_{i_1},{\cdots},{\rm link}_{L_{i_s}}\!\sigma_{i_s})*
\Delta\!^{[n']},\,
\{i_1,{\cdots},i_s\}=[m]\setminus\hat\sigma.$$
This implies that
${\cal S}(K;L_1,{\cdots},L_m)$ is a homological sphere ($L$ is a homological sphere if ${\rm link}_L\sigma$ is homology spherical
for all $\sigma\in L$) if and only if $K$ is a homological sphere and each $L_k=\partial\Delta\!^{[n_k]}$.
This result is a part of Theorem~6.6 of \cite{AY}.
\vspace{3mm}

We have ring isomorphism $H^*({\cal S}(K;\underline{X},\underline{A}))
\cong H^*(|{\cal S}(K;\underline{X},\underline{A})|)$, where $|\cdot|$ means geometrical
realization.
So $\w H^{*}({\cal S}(K;\underline{X},\underline{A}))$ is
a ring by adding a unit to it.
This ring is not considered in this paper. \vspace{3mm}

{\bf Theorem 3.11}{\it\, For a total homology split ${\cal S}(K;\underline{X},\underline{A})$, we have
$$\begin{array}{l}
\quad H_*^{\SS_n}({\cal S}(K;\underline{X},\underline{A}))\vspace{2mm}\\
\cong H_*^{\XX_m}(K)\,\widehat\otimes\,
\big(H_*^{\XX_m;\SS_n}(\underline{X},\underline{A})\big)\vspace{2mm}\\
\cong H_*^{\XX_m}(K)\,\widehat\otimes\,
\big(H_*^{\XX\!;\SS_{n_1}}(X_1,A_1){\otimes}{\cdots}{\otimes}H_*^{\XX\!;\SS_{n_m}}(X_m,A_m)\big),
\end{array}$$
$$\begin{array}{l}
\quad H^{\,*}_{\!\SS_n}({\cal S}(K;\underline{X},\underline{A}))\vspace{2mm}\\
\cong H^{\,*}_{\!\XX_m}(K)\,\widehat\otimes\,
\big(H^{\,*}_{\!\XX_m;\SS_n}(\underline{X},\underline{A})\big)\vspace{2mm}\\
\cong H^{\,*}_{\!\XX_m}(K)\,\widehat\otimes\,
\big(H^{\,*}_{\!\XX\!;\SS_{n_1}}(X_1,A_1){\otimes}{\cdots}{\otimes}H^{\,*}_{\!\XX\!;\SS_{n_m}}(X_m,A_m)\big),  \end{array}
$$
where $\SS_n=\SS_{n_1}{\times}{\cdots}{\times}\SS_{n_m}$ by the 1-1 correspondence
$$(\sigma\!,\,\omega)\leftrightarrow (\sigma_1,\omega_1,{\cdots},\sigma_m,\omega_m),\,\sigma_k=\sigma{\cap}[n_k],\,\omega_k=\omega{\cap}[n_k].$$
$\SS=\XX$ or $\RR$ ($\SS_{n_1}=\XX_{n_1}$, $\SS_{n_2}=\RR_{n_2}$ is possible).
Since $\XX_m=\XX{\times}{\cdots}{\times}\XX$ as in Definition 3.5, we have
$(\XX_m;\SS_n)=(\XX;\SS_{n_1}){\times}{\cdots}{\times}(\XX;\SS_{n_m})$ and so for $(\hat\sigma,\hat\omega)\in\XX_m$, $(\sigma\!,\,\omega)\in\SS_n$,
we have by definition
$$H_*^{\Lambda;\Gamma}(-)=\oplus_{\alpha\in\Lambda,\,\beta\in\Gamma}\,H_*^{\alpha;\beta}(-),\quad
H^*_{\Lambda;\Gamma}(-)=\oplus_{\alpha\in\Lambda,\,\beta\in\Gamma}\,H^*_{\alpha;\beta}(-),$$
$$H_*^{s;\sigma_k,\omega_k}(X_k,A_k)=
\left\{\begin{array}{cl}
\Sigma{\rm coker}\,\iota_{\sigma_k,\omega_k}&{\rm if}\, s={\scriptstyle(\{1\},\emptyset)},\vspace{1mm}\\
\Sigma{\rm ker}\,\iota_{\sigma_k,\omega_k}&{\rm if}\, s={\scriptstyle(\emptyset,\{1\})},\vspace{1mm}\\
\Sigma{\rm im}\,\iota_{\sigma_k,\omega_k}&{\rm if}\, s={\scriptstyle(\,\emptyset,\,\emptyset)},
\end{array}\right.$$
$$H^*_{s;\sigma_k,\omega_k}(X_k,A_k)=
\left\{\begin{array}{cl}
\Sigma{\rm ker}\,\iota_{\sigma_k,\omega_k}^\circ&{\rm if}\, s={\scriptstyle(\{1\},\emptyset)},\vspace{1mm}\\
\Sigma{\rm coker}\,\iota_{\sigma_k,\omega_k}^\circ&{\rm if}\, s={\scriptstyle(\emptyset,\{1\})},\vspace{1mm}\\
\Sigma{\rm im}\,\iota_{\sigma_k,\omega_k}^\circ&{\rm if}\, s={\scriptstyle(\,\emptyset,\,\emptyset)},
\end{array}\right.$$
$$H_*^{\hat\sigma\!,\,\hat\omega;\sigma\!,\,\omega}(\underline{X},\underline{A})
=H_1{\otimes}{\cdots}{\otimes}H_m,\quad
H_k=\left\{\begin{array}{rl}
\Sigma{\rm coker}\,\iota_{\sigma_k,\omega_k}&{\rm if}\, k\in\hat\sigma,\vspace{1mm}\\
\Sigma{\rm ker}\,\iota_{\sigma_k,\omega_k}&{\rm if}\, k\in\hat\omega,\vspace{1mm}\\
\Sigma{\rm im}\,\iota_{\sigma_k,\omega_k}
&{\rm otherwise},
\end{array}\right.$$
$$H^*_{\hat\sigma\!,\,\hat\omega;\sigma\!,\,\omega}(\underline{X},\underline{A})
=H^1{\otimes}{\cdots}{\otimes}H^m,\quad
H^k=\left\{\begin{array}{rl}
\Sigma{\rm ker}\,\iota_{\sigma_k,\omega_k}^\circ&{\rm if}\, k\in\hat\sigma,\vspace{1mm}\\
\Sigma{\rm coker}\,\iota_{\sigma_k,\omega_k}^\circ&{\rm if}\, k\in\hat\omega,\vspace{1mm}\\
\Sigma{\rm im}\,\iota_{\sigma_k,\omega_k}^\circ&{\rm otherwise},
\end{array}\right.$$
where $\iota_{\sigma_k,\omega_k}$ is as in Definition 3.8
with dual $\iota_{\sigma_k,\omega_k}^\circ$ and $\Sigma$ means suspension.

If each $\iota_{\sigma_k,\omega_k}$ is an epimorphism, then all $\XX$ is replaced by $\RR$ and we have
$H_*^{\RR;\SS_{n_k}}(X_k,A_k)=H_*^{\SS_{n_k}}(A_k)$,
$H^{\,*}_{\RR;\SS_{n_k}}(X_k,A_k)=H^{\,*}_{\SS_{n_k}}(A_k)$.

If the (right) total (co)homology group is taken over a field, then the theorem holds for all
polyhedral join complexes.
\vspace{2mm}

Proof}\, A corollary of Theorem 3.7
by taking $\underline{\vartheta}=\{\vartheta_k\}_{k=1}^m$ with split inclusion
$\vartheta_k\colon T_*^{\SS_{n_k}}(A_k)\to T_*^{\SS_{n_k}}(X_k)$
the (right) total chain complex inclusion.
Then $(C_*(K;\underline{\vartheta}),d)=(T_*^{\SS_n}({\cal S}(K,\underline{X},\underline{A})),d)$
and $H_*^\XX\!(\theta_k)=H_*^{\XX\!;\SS_k}(\theta_k)=H_*^{\XX\!;\SS_k}(X_k,A_k)$.
\hfill$\Box$\vspace{3mm}

{\bf Example 3.12}  Apply Theorem 3.11 for ${\cal S}(K;\underline{X},\underline{A})={\cal S}(K;L_1,{\cdots},L_m)$,
$L_k\neq\{\,\}$ for $k=1,{\cdots},m$.
So either all $H_*^{\RR_{n_k}}(L_k)$ are free or the (co)homology is taken over a field.

For $\SS=\RR$, $\theta_k\colon H_*^{\RR_{n_k}}(L_k)\to H_*^{\RR_{n_k}}(\Delta\!^{[n_k]})\cong \Bbb Z$
is an epimorphism. So
$$H_*^{\RR_n}({\cal S}(K;L_1,{\cdots},L_m))\cong H_*^{\RR_m}(K)\,\widehat\otimes\,
\big(H_*^{\RR_{n_1}}(L_1){\otimes}{\cdots}{\otimes}H_*^{\RR_{n_m}}(L_m)\big),$$
$$H^*_{\RR_n}({\cal S}(K;L_1,{\cdots},L_m))\cong H^*_{\RR_m}(K)\,\widehat\otimes\,
\big(H^*_{\RR_{n_1}}(L_1){\otimes}{\cdots}{\otimes}H^*_{\RR_{n_m}}(L_m)\big).$$

For $\SS=\XX$, $\theta_k\colon H_*^{\XX_{n_k}}\!(L_k)\to H_*^{\XX_{n_k}}\!(\Delta\!^{[n_k]})$
is not an epimorphism.
By definition, $H_*^{\emptyset,\{1\};\sigma_k,\omega_k}(\Delta\!^{[n_k]}\!{,}L_k)\cong H_*^{\sigma_k,\omega_k}(L_k)$ for $\omega_k\neq\emptyset$,
$H_*^{\emptyset,\emptyset;\sigma_k,\emptyset}(\Delta\!^{[n_k]}\!{,}L_k)\cong\Bbb Z$ for $\sigma_k\in L_k$ and
$H_*^{\{1\},\emptyset;\sigma_k,\emptyset}(\Delta\!^{[n_k]}\!{,}L_k)\cong\Bbb Z$ for $\sigma_k\notin L_k$. So
by identifying $X{\otimes}\Bbb Z{\otimes}Y$ with  $X{\otimes}Y$, we have that for $(\sigma\!,\,\omega)\in\XX_n$,
$$H_*^{\sigma\!,\,\omega}({\cal S}(K;L_1,{\cdots},L_m))\cong H_*^{\hat\sigma\!,\,\hat\omega}(K)\otimes
\big(\otimes_{\omega_k\neq\emptyset}\,H_*^{\sigma_k,\omega_k}(L_k)\big),$$
$$H^*_{\sigma\!,\,\omega}({\cal S}(K;L_1,{\cdots},L_m))\cong H^*_{\hat\sigma\!,\,\hat\omega}(K)\otimes
\big(\otimes_{\omega_k\neq\emptyset}\,H^*_{\sigma_k,\omega_k}(L_k)\big),$$
where $\sigma_k=\sigma{\cap}[n_k]$, $\omega_k=\omega{\cap}[n_k]$, $\hat\sigma=\{k\,|\,\sigma_k\!\notin\! L_k,\,\omega_k=\emptyset\}$, $\hat\omega=\{k\,|\,\omega_k\!\neq\!\emptyset\}$.

\section{Universal Algebra}

In this section, we prove that the cohomology algebra of
$C_*(K;\underline{\vartheta})$ is the diagonal tensor product of the total cohomology algebra of
$K$ induced by $\underline{\vartheta}$ and the character algebra of the dual $\underline{\vartheta}^\circ$
in Theorem 4.7.
The (co)associativity is not required for a (co)algebra as in \cite{Z} .
\vspace{3mm}

{\bf Theorem 4.1}\, {\it Let $\vartheta\colon (C_*,d)\to(D_*,d)$ be a split inclusion
with induced homology homomorphism $\theta\colon U_*\to V_*$
such that $\vartheta$ is also a coalgebra homomorphism $\vartheta\colon (C_*,\psi_C)\to(D_*,\psi_D)$
with induced homology coalgebra homomorphism $\theta\colon (U_*,\psi_U)\to(V_*,\psi_V)$.

Then the group $C_*^{\XX}\!(\theta)$ in Definition 3.3 has a unique character coproduct
$\widehat\psi(\vartheta)$ satisfying the following three conditions.

i) $\widehat\psi(\vartheta)$ makes the following diagram ($q$, $q'$ and $\vartheta'$ as in Theorem 3.4).

\begin{center}
\setlength{\unitlength}{1mm}
\begin{picture}(79,37)
\put(0,20){${\scriptstyle C_*}$}\put(30,20){${\scriptstyle U_*}$}
\put(10,30){${\scriptstyle C_*{\otimes}C_*}$}\put(42,30){${\scriptstyle U_*{\otimes}U_*}$}
\put(0,0){${\scriptstyle D_*}$}\put(30,0){${\scriptstyle C_*^{\XX\!}}$}
\put(10,10){${\scriptstyle D_*{\otimes}D_*}$}\put(41,10){${\scriptstyle C_*^{\XX\!}{\otimes}C_*^{\XX\!}}$}
\thicklines
\put(6,21){\vector(1,0){22}}\put(12,20.5){$^{q'}$} \put(2,18){\vector(0,-1){13}}\put(2.5,12){$^\vartheta$}
\put(6,1){\vector(1,0){22}}\put(12,0.5){$^{q}$} \put(32,18){\vector(0,-1){13}}\put(32.5,12){$^{\vartheta'}$}
\thinlines
\put(23,31){\vector(1,0){15}}\put(24,30.5){$^{q'{\otimes}q'}$}
\put(15.5,28){\line(0,-1){6}}\put(15.5,20){\vector(0,-1){6}}\put(16,22){$^{\vartheta{\otimes}\vartheta}$}
\put(23,11){\line(1,0){8}}\put(33,11){\vector(1,0){6}}\put(24,10.5){$^{q{\otimes}q}$}
\put(47.5,28){\vector(0,-1){14}}\put(48,22){$^{\vartheta'{\otimes}\vartheta'}$}
\thicklines
\put(4,24){\vector(1,1){5}}\put(7,23.5){$^{\psi_C}$} \put(35,24){\vector(1,1){5}}\put(38,23){$^{\psi_U}$}
\put(4,3.5){\vector(1,1){5}}\put(7,2.5){$^{\psi_D}$} \put(35,3.5){\vector(1,1){5}}\put(38,2){$^{\widehat\psi(\vartheta)}$}
\put(57,0.5){${\scriptstyle (C_*^\XX=\,C_*^{\XX}\!(\theta))}$}
\end{picture}\end{center}
commutative except the homotopy commutative $(q{\otimes}q)\psi_D\simeq \widehat\psi(\vartheta)q$.

ii) $\widehat\psi(\vartheta)$ is independent of the choice of $\psi_C,\psi_D$ up to homotopy,
i.e., if $\psi_C,\psi_D$ are replaced by $\psi'_C,\psi'_D$ such that $\psi'_C\simeq\psi_C$, $\psi'_D\simeq\psi_D$
and we get $\widehat\psi'(\vartheta)$ for $\psi'_C$ and $\psi'_D$,
then $\widehat\psi'(\vartheta)=\widehat\psi(\vartheta)$.

iii) Denote by $\alpha={\rm coker}\,\theta$, $\beta=\Sigma\,{\rm ker}\,\theta$,
$\gamma={\rm ker}\,\theta$, $\eta={\rm im}\,\theta$.
Then $\widehat\psi(\vartheta)$ satisfies the following four conditions.

(1) $\widehat\psi(\vartheta)(\eta)\subset \eta{\otimes}\eta\oplus\gamma{\otimes}\eta\oplus \eta{\otimes}\gamma
\oplus \gamma{\otimes}\gamma$.

(2) $\widehat\psi(\vartheta)(\gamma)\subset \gamma{\otimes}\gamma\oplus \gamma{\otimes}\eta\oplus \eta{\otimes}\gamma$.

(3) $\widehat\psi(\vartheta)(\beta)\subset\big(\beta{\otimes}\gamma\oplus \beta{\otimes}\eta\oplus\eta{\otimes}\beta\big)
\oplus\big(\alpha{\otimes}\alpha\oplus \alpha{\otimes}\eta\oplus
\eta{\otimes}\alpha\oplus \eta{\otimes}\eta\big)$.

(4) $\widehat\psi(\vartheta)(\alpha)\subset \alpha{\otimes}\alpha\oplus \alpha{\otimes}\eta\oplus \eta{\otimes}\alpha
\oplus \eta{\otimes}\eta$.
\vspace{2mm}

Proof}\, Similar to that of Theorem 2.8 of \cite{Z}.
\hfill$\Box$\vspace{3mm}

{\bf Definition 4.2} For the $\vartheta$ in Theorem 4.1, all the chain complexes in Theorem 3.4 are coalgebras defined as follows.

The {\it character coalgebra complex} of $\vartheta$ is $(C_*^\XX\!(\theta),\widehat\psi(\vartheta))$.

If $\theta$ is an epimorphism, then $C_*^\XX(\theta)=C_*^\RR(\theta)$ and $(C_*^\RR(\theta),\widehat\psi(\vartheta))$
is called the {\it right character coalgebra complex} of $\vartheta$.

The {\it character coalgebra} of $\vartheta$ is $(H_*^{\XX}\!(\theta),\psi(\vartheta))$ with coproduct
defined as follows.

(1) $\psi(\vartheta)(x)=\widehat\psi(\vartheta)(x)=\psi_V(x)$ for all $x\in \alpha={\rm  coker}\,\theta$.

(2) $\psi(\vartheta)(y)=\widehat\psi(\vartheta)(y)=\psi_U(y)$ for all $y\in\eta={\rm coim}\,\theta$.

(3) For $z\in\gamma={\rm ker}\,\theta$, denote by $\overline z$ be the unique element in $\beta$ such that $d\,\overline z=z$.
Suppose $\hat\psi(\vartheta)(z)=\psi_U(z)=\Sigma y_i{\otimes}z_i+\Sigma z'_j{\otimes} z''_j$
with $z_i,z'_j,z''_j\in\gamma$ and $y_i\in\eta$,
then $\hat\psi(\vartheta)(\overline z)=
\Sigma{\scriptstyle(-1)}^{|y_i|}y_i{\otimes}\overline z_i+\Sigma\overline z'_j{\otimes} z''_j+\zeta$
with $\zeta\in(\alpha{\oplus}\eta){\otimes}(\alpha{\oplus}\eta)$.
Define $\psi(\vartheta)(z)=\widehat\psi(\vartheta)(z)+\zeta$.

If $\theta$ is an epimorphism, then $(H_*^\RR(\theta),\psi(\vartheta))$ is the {\it right character coalgebra} of $\vartheta$.
The group isomorphisms $H_*^\RR(\theta)\cong V_*$ is in general not an algebra isomorphism.

For $\SS=\XX$ or $\SS$, the dual algebra $(H^{\,*}_{\SS}(\theta^\circ),\pi(\vartheta^\circ))$)
of $(H_*^\SS(\theta),\psi(\vartheta))$ is the {\it (right) character algebra} of $\vartheta^\circ$.

The {\it index coalgebra complex}
$(T_*^\XX\!,\widehat\psi_\vartheta)$ of $\vartheta$ is defined as follows.
Let symbols $x,x'_1,x''_2,\cdots$ be one of $\alpha,\beta,\gamma,\eta$
that are both the generators of $T_*^\XX$ and the group summands of $C_*^\XX\!(\theta)$ in
Theorem 4.1. If the group summand $x$ satisfies $\widehat\psi(\vartheta)(x)\subset\oplus_i\,x'_i{\otimes}x''_i$
such that no summand $x'_i{\otimes}x''_i$ can be canceled,
then the generator $x$ satisfies $\widehat\psi_\vartheta(x)=\Sigma_i\,x'_i{\otimes}x''_i$.

If $\theta$ is an epimorphism, then $(T_*^\RR,\widehat\psi_{\vartheta})$ is
the {\it right index coalgebra complex} of $\vartheta$.
\vspace{3mm}

Notice that $\hat\psi(\vartheta)$ is degree preserving but $\psi(\vartheta)$ is in general not.
When $\psi(\vartheta)=\hat\psi(\vartheta)$
and so $\psi(\vartheta)$ keeps degree, the coproduct is called normal,
for example, (right) normal, (right) strictly normal.
\vspace{3mm}

{\bf Example 4.3} Take $\vartheta\colon(T_*^{\RR_m}(L),\psi_C)\to(T_*^{\RR_m}(\Delta^{[m]}),\psi_D)$ in Theorem~4.1 to be the
right total chain complex homomorphism induced by inclusion such that $T_*^{\RR_m}(L)$ is free and the vertex set of
$L$ is $[m]$. $\psi_C,\psi_D$ are the right universal coproduct defined in Definition~7.4 of \cite{Z}.
Then the induced homology coalgebra homomorphism
$\theta\colon H_*^{\RR_m}(L)\to H_*^{\RR_m}(\Delta^{[m]})\cong\Bbb Z$
is an epimorphism. By definition,
$\eta\cong H_*^{\emptyset,\emptyset}(L)\cong\Bbb Z$ with generator denoted by $1$,
$\gamma\cong \oplus_{\omega\neq\emptyset}H_*^{\emptyset,\omega}(L)$, $\alpha=0$.
Since there is no ghost vertex, $H_0^{\emptyset,\omega}(L)=0$ for all $\omega\neq\emptyset$.
So $\widehat\psi(\vartheta)(1)=1{\otimes}1$. For $x\in\gamma$,
$\widehat\psi(\vartheta)(x)=1{\otimes}x+x{\otimes}1+\Sigma x'_i{\otimes}x''_i$ with $x'_i,x''_i\in\gamma$.
Since the degree of $1$ is $0$, $1{\otimes}1$ can not be a summand of $\widehat\psi(\vartheta)(\overline x)$.
So $\widehat\psi(\vartheta)(\overline x)=1{\otimes}\overline x+\overline x{\otimes}1+\Sigma\overline x'_i{\otimes}x''_i$.
With these, we have

$\widehat\psi(\vartheta)(\eta)\subset\eta{\otimes}\eta$ ($\eta={\rm im}\theta$) and so
$\widehat\psi_{\vartheta}(\eta)=\eta{\otimes}\eta$ ($\eta$ the generator of $T_*^{\RR}$).

$\widehat\psi(\vartheta)(\gamma)\subset\gamma{\otimes}\gamma\oplus\gamma{\otimes}\eta\oplus\eta{\otimes}\gamma$ and so
$\widehat\psi_{\vartheta}(\gamma)=\gamma{\otimes}\gamma+\gamma{\otimes}\eta+\eta{\otimes}\gamma$.

$\widehat\psi(\vartheta)(\beta)\subset\beta{\otimes}\gamma\oplus\beta{\otimes}\eta\oplus\eta{\otimes}\beta$ and so
$\widehat\psi_{\vartheta}(\beta)=\beta{\otimes}\gamma+\beta{\otimes}\eta+\eta{\otimes}\beta$.

So  $\widehat\psi_{\vartheta}$ is the right strictly normal coproduct in Definition 7.4 of \cite{Z}
and we have the coalgebra isomorphism $(H_*^\RR(\theta),\psi(\vartheta))=(H_*^{\RR_m}(L),\psi_U)$.
\vspace{3mm}

{\bf Example 4.4} Take $\vartheta\colon(T_*^{\XX_m}(L),\psi_C)\to(T_*^{\XX_m}(\Delta\!^{[m]}),\psi_D)$ in Theorem~4.1 to be the
total chain complex homomorphism induced by inclusion  such that $T_*^{\XX_m}(L)$ is free and the vertex
set of $L$ is $[m]$. $\psi_C,\psi_D$ are the universal coproduct defined in Definition~6.1 of \cite{Z}.
The induced homology coalgebra homomorphism $\theta\colon H_*^{\XX_m}(L)\to H_*^{\XX_m}(\Delta\!^{[m]})$
is not an epimorphism. By definition,
$\eta\cong\oplus_{\sigma\in K} H_*^{\sigma,\emptyset}(L)$, where each $H_*^{\sigma,\emptyset}(L)\cong\Bbb Z$
with generator denoted by $\eta_\sigma$.
$\alpha\cong\oplus_{\sigma\notin K} H_*^{\sigma,\emptyset}(\Delta\!^{[n_k]})$, where each $H_*^{\sigma,\emptyset}(\Delta\!^{[n_k]})\cong\Bbb Z$
with generator denoted by $\alpha_\sigma$.
$\gamma\cong \oplus_{\omega\neq\emptyset}\,H_*^{\sigma\!,\,\omega}(L)$.
Since there is no ghost vertex, $H_0^{\sigma\!,\,\omega}(L)=0$ for all $\omega\neq\emptyset$. So
$\widehat\psi(\vartheta)(\alpha_\sigma)=\Sigma(\eta_{\sigma'}{\otimes}\eta_{\sigma''})
+(\alpha_{\sigma'}{\otimes}\eta_{\sigma''})+(\eta_{\sigma'}{\otimes}\alpha_{\sigma''})+(\alpha_{\sigma'}{\otimes}\alpha_{\sigma''})$,
$\widehat\psi(\vartheta)(\eta_\sigma)=\Sigma(\eta_{\sigma'}{\otimes}\eta_{\sigma''})$,
where the sums are taken over all $\sigma'{\cup}\sigma''\subset\sigma$.
For $x\in\gamma$,
$\widehat\psi(\vartheta)(x)=\Sigma_{\sigma'\subset\sigma}(\eta_{\sigma'}{\otimes}x{+}x{\otimes}\eta_{\sigma'})+\Sigma x'_i{\otimes}x''_i$ with $x'_i,x''_i\in\gamma$.
Since the degree of all $\eta_\sigma$ and $\alpha_\sigma$ is $0$, $\xi_{\sigma'}{\otimes}\xi_{\sigma''}$ ($\xi=\alpha$ or $\eta$)
can not be a summand of $\widehat\psi(\vartheta)(\overline x)$.
So $\widehat\psi(\vartheta)(\overline x)=\Sigma_{\sigma'\subset\sigma}(\eta_{\sigma'}{\otimes}\overline x{+}\overline x{\otimes}\eta_{\sigma'})
+\Sigma\overline x'_i{\otimes}x''_i$.
With these, we have

$\widehat\psi(\vartheta)(\eta)\subset\eta{\otimes}\eta$ and so
$\widehat\psi_{\vartheta}(\eta)=\eta{\otimes}\eta$.

$\widehat\psi(\vartheta)(\gamma)\subset\gamma{\otimes}\gamma\oplus\gamma{\otimes}\eta\oplus\eta{\otimes}\gamma$ and so
$\widehat\psi_{\vartheta}(\gamma)=\gamma{\otimes}\gamma+\gamma{\otimes}\eta+\eta{\otimes}\gamma$.

$\widehat\psi(\vartheta)(\beta)\subset\beta{\otimes}\gamma\oplus\beta{\otimes}\eta\oplus\eta{\otimes}\beta$ and so
$\widehat\psi_{\vartheta}(\beta)=\beta{\otimes}\gamma+\beta{\otimes}\eta+\eta{\otimes}\beta$.

$\widehat\psi(\vartheta)(\alpha)\subset(\alpha{\oplus}\eta){\otimes}(\alpha{\oplus}\eta)$ and so
$\widehat\psi_{\vartheta}(\alpha)=\alpha{\otimes}\alpha+\alpha{\otimes}\eta+\eta{\otimes}\alpha+\eta{\otimes}\eta$.

We call the above $\widehat\psi_{\vartheta}$ {\it the strictly normal coproduct} of $T_*^\XX$.
The coproduct $\psi(\vartheta)=\widehat\psi(\vartheta)$.
$H_*^{\XX_m}(L)$ is a subcoalgebra of
$H_*^{\XX}(\theta)=\alpha\oplus H_*^{\XX_m}(L)$. Dually,
we call $(H^*_{\XX}(\theta^\circ),\pi(\vartheta^\circ))$ the {\it augmented cohomology algebra} of
$L$ relative to $[m]$ and denote it by $\widetilde H^*_{\XX_m}(L)$.
Then $\widetilde H^*_{\XX}(L)=\alpha^\circ\oplus H^*_{\XX_m}(L)$ with $\alpha^\circ$
the dual group of $\alpha$ and $\alpha^\circ$ is an ideal such that
$H^*_{\XX_m}(L)=\widetilde H^*_{\XX_m}(L)/\alpha^\circ$.
\vspace{3mm}

Notice that coalgebras $(T_*^\SS,\hat\psi_{\vartheta})$ and
$(H_*(\theta),\psi(\vartheta))$ in Example 3.3 and 3.4 remain unchanged
if the (right) universal coproducts $\psi_C,\psi_D$ are replaced by any other coproduct.
\vspace{3mm}

{\bf Theorem 4.5} {\it For the coalgebras in Definition 4.2,
all the chain homomorphisms in Theorem 3.4 induce cohomology algebra isomorphisms.
\vspace{2mm}

Proof}\, Similar to that of Theorem 6.4 of \cite{Z}.
\hfill$\Box$\vspace{3mm}

{\bf Definition 4.6} Let $K$, $\underline{\vartheta}$, $\underline{\theta}$ be
as in Definition 3.5 and Definition 3.6
such that each $\vartheta_k\colon ((C_k)_*,\psi_{C_k})\to((D_k)_*,\psi_{D_k})$
and $\theta_k\colon ((U_k)_*,\psi_{U_k})\to((V_k)_*,\psi_{V_k})$
satisfy the condition of Theorem 4.1. $\SS=\XX$ or $\RR$.
Then all the chain complexes in Theorem 3.7 are coalgebras defined as follows.

The polyhedral product inclusion complex $C_*(K;\underline{\vartheta})$ is a subcoalgebra of $(D_1)_*{\otimes}{\cdots}{\otimes}(D_m)_*$.
Its cohomology algebra product is denoted by
$\cup_{(K;\underline{\vartheta}^\circ)}$.

The (right) polyhedral product character complex $C_*^{\SS_m}(K;\underline{\vartheta})$ is a subcoalgebra of
the {\it (right) character coalgebra complex} of $\underline{\vartheta}$, which is denoted by
$(C_*^{\SS_m}(\underline{\vartheta}),\widehat\psi(\underline{\vartheta}))=
(C_*^{\SS}\!(\vartheta_1){\otimes}{\cdots}{\otimes} C_*^{\SS}\!(\vartheta_m),\widehat\psi(\vartheta_1){\otimes}{\cdots}{\otimes}
\widehat\psi(\vartheta_m))$.

The {\it (right) character coalgebra} of $\underline{\vartheta}$ is
$(H_*^{\SS_m}(\underline{\theta}),\psi(\underline{\vartheta}))=
(H_*^{\SS}\!(\theta_1)\otimes{\cdots}\otimes H_*^{\SS}\!(\theta_m),\psi(\vartheta_1){\otimes}{\cdots}{\otimes}
\psi(\vartheta_m))$. The {\it (right) character algebra} of $\underline{\vartheta}^\circ$ is the dual algebra
$(H^*_{\SS_m}(\underline{\theta}^\circ),\pi(\underline{\vartheta}^\circ))=
(H^*_{\SS}(\theta_1^\circ)\otimes{\cdots}\otimes H^*_{\SS}(\theta_m^\circ),\pi(\vartheta_1^\circ){\otimes}{\cdots}{\otimes}
\pi(\vartheta_m^\circ))$.

The (right) total chain complex $T_*^{\SS_m}(K)$ of $K$ is a subcoalgebra
the {\it (right) index coalgebra complex} $(T_*^{\SS_m}{,}\widehat\psi_{\underline{\vartheta}})=
(T_*^{\SS}\!{\otimes}{\cdots}{\otimes}T_*^{\SS},\widehat\psi_{\vartheta_1}\!{\otimes}{\cdots}{\otimes}
\widehat\psi_{\vartheta_m})$ of $\underline{\vartheta}$
and is called the {\it (right) total coalgebra complex} of $K$ induced by $\underline{\vartheta}$.
Its cohomology algebra is called the {\it (right) total cohomology algebra} of $K$ induced by $\underline{\vartheta}$ and
is denoted by $(H^*_{\SS_m}(K),\cup_{(K;\underline{\vartheta}^\circ)})$.
\vspace{3mm}

{\bf Theorem 4.7} {\it For the coalgebras in Definition 4.6,
all the chain homomorphisms in Theorem 3.7 induce cohomology algebra isomorphisms.
So we have cohomology algebra isomorphism
$$(H^*(C^*(K;\underline{\theta}^\circ)),\cup_{(K;\underline{\vartheta}^\circ)})
\cong (H^*_{\XX_m}(K)\,\widehat\otimes\, H^*_{\XX_m}(\underline{\theta}^\circ),\cup_{(K;\underline{\vartheta}^\circ)}\,\widehat\otimes\,
\pi(\underline{\vartheta}^\circ)).$$

Proof}\, The group isomorphisms in Theorem 3.7 naturally induce cohomology algebra isomorhisms.
\hfill$\Box$\vspace{3mm}

{\bf Example 4.8} Let everything be as in Theorem 2.9.
We compute the cohomology algebra over a field of ${\cal Z}(K;\underline{Y},\underline{B})$.
By Theorem 6.9 and 7.11 of \cite{Z}, we have algebra isomorphisms
$$\begin{array}{c}
(H^*(Y_k),\cup)\cong\big(H^*_{\SS_{n_k}}(X_k)\,\widehat\otimes\,H^*_{\SS_{n_k}}(\underline{U_k},\underline{C_k}),
\cup_{X_k}\,\widehat\otimes\,\pi_{(\underline{U_k},\underline{C_k})}\big),\vspace{2mm}\\
(H^*(B_k),\cup)\cong\big(H^*_{\SS_{n_k}}(A_k)\,\widehat\otimes\,H^*_{\SS_{n_k}}(\underline{U_k},\underline{C_k}),
\cup_{A_k}\,\widehat\otimes\,\pi_{(\underline{U_k},\underline{C_k})}\big),
\end{array}$$
where $\SS_{n_k}=\XX_{n_k}$ or $\RR_{n_k}$,
$\cup_{X_k}$ and $\cup_{A_k}$ are the (right) universal (or (right) normal, etc.) product appearing in the theorems
induced by the coproduct $\psi_{X_k}$ and $\psi_{A_k}$.
Take $\vartheta_k\colon(T_*^{\SS_{n_k}}(A_k),\psi_{A_k})\to(T_*^{\SS_{n_k}}(X_k),\psi_{X_k})$ to be the (right) total
chain complex inclusion and apply Theorem 4.5 for $\underline{\vartheta}=\{\vartheta_k\}_{k=1}^m$.
Then we have algebra isomorphisms
$$\,\,(H^*({\cal Z}(K;\underline{Y},\underline{B}),\cup)\cong\big(H^*_{\SS_n}({\cal S}(K;\underline{X},\underline{A}))\,\widehat\otimes\,
H^*_{\SS_n}(\underline{U},\underline{C}),
\cup\,\widehat\otimes\,\pi_{(\underline{U},\underline{C})}\big),$$
$$(H^*_{\SS_n}({\cal S}(K;\underline{X},\underline{A})),\cup)\cong
\big(H^*_{\XX_m}(K)\,\widehat\otimes\,H^*_{\XX_m;\SS_n}(\underline{X},\underline{A}),
\cup_{(K;\underline{\vartheta}^\circ)}\,\widehat\otimes\,\pi(\underline{\vartheta}^\circ)\big),$$
where $\SS_n=\SS_{n_1}{\times}{\cdots}{\times}\SS_{n_m}$ ($\SS_{n_1}=\XX_{n_1}$, $\SS_{n_2}=\RR_{n_2}$ is possible), i.e.,
$$H^*_{\SS_n}({\cal Z}(K;\underline{Y},\underline{B}))\cong
H^*_{\XX_m}(K)\,\widehat\otimes\,H^*_{\XX_m;\SS_n}(\underline{X},\underline{A})
\,\widehat\otimes\,H^*_{\SS_n}(\underline{U},\underline{C}).$$

When all $\SS_{n_k}=\XX_{n_k}$ and the coproduct $\psi_{X_k},\psi_{A_k}$ are
the universal (or normal, etc.) coproduct, then
$(H^*_{\XX_n}({\cal S}(K;\underline{X},\underline{A})),\cup_{(K;\underline{\vartheta}^\circ)})$ is
just the universal (or normal, etc.) algebra of ${\cal S}(K;\underline{X},\underline{A})$.
The right algebra case also holds by replacing $\XX$ by $\RR$.
\vspace{3mm}

{\bf Example 4.9} Let $L_k$ satisfy the condition of Example 4.3 and 4.4 for $k=1,{\cdots},m$.
Then we have algebra isomorphisms
$$H^*_{\RR_n}({\cal S}(K;L_1,{\cdots},L_m))
\cong H^*_{\RR_m}(K){\widehat\otimes}
\big(H^*_{\RR_{n_1}}(L_1){\otimes}{\cdots}{\otimes}H^*_{\RR_{n_m}}(L_m)\big),$$
$$H^*_{\XX_n}({\cal S}(K;L_1,{\cdots},L_m))
\cong H^*_{\XX_m}(K){\widehat\otimes}
\big(\widetilde H^*_{\XX_{n_1}}(L_1){\otimes}{\cdots}{\otimes}\widetilde H^*_{\XX_{n_m}}(L_m)\big),$$
where $H^*_{\XX_m}(K)$ ($H^*_{\RR_m}(K)$) is the {\it (right) strictly normal algebra} of $K$
induced by the (right) strictly normal coproduct of $T_*^\XX$ ($T_*^\RR$)
in Example 4.4 (4.3) by Definition 4.6, $\widetilde H^*_{\XX_{n_k}}(L_k)$ is the augmented cohomology algebra of $L_k$
relative to $[n_k]$ defined in Example 4.4.
\vspace{3mm}

\section{Duality Isomorphism}\vspace{3mm}

In this section, we compute the Alexander duality isomorphism
on some special type of polyhedral product spaces.\vspace{3mm}

{\bf Theorem 5.1} {\it Let $(\underline{X},\underline{A})=\{(X_k,A_k)\}_{k=1}^m$ be a sequence of topological
pairs satisfying the following conditions.

1) Each homology group homomorphism $i_k\colon H_*(A_k)\to H_*(X_k)$
induced by inclusion is a split homomorphism.

2) Each $X_k$ is a closed orientable manifold of dimension $r_k$.

3) Each $A_k$ is a proper compact polyhedron subspace of $X_k$.

Let $(\underline{X},\underline{A}^c)=\{(X_k,A_k^c)\}_{k=1}^m$ with $A_k^c=X_k{\setminus}A_k$.
Then for all $(\sigma\!,\omega)\in\XX_m$, there are duality isomorphisms (\,$r=r_1{+}{\cdots}{+}r_m$\,)
$$\gamma_{\sigma\!,\,\omega}\colon H_*^{\sigma\!,\,\omega}(\underline{X},\underline{A})\to H^{r-|\omega|-*}_{\tilde\sigma\!,\,\omega}(\underline{X},\underline{A}^c),$$
$$\gamma^\circ_{\sigma\!,\,\omega}\colon H^*_{\sigma\!,\,\omega}(\underline{X},\underline{A})\to H_{r-|\omega|-*}^{\tilde\sigma\!,\,\omega}(\underline{X},\underline{A}^c),$$
where $\w\sigma=[m]{\setminus}(\sigma{\cup}\omega)$,
$H_*^{\sigma\!,\,\omega}(-)$ and $H^*_{\sigma\!,\,\omega}(-)$ are as in Theorem 3.9.

If the (co)homology is taken over a field, then the conclusion holds for
$(\underline{X},\underline{A})$ satisfying the following conditions.

1) Each $X_k$ is a closed manifold of dimension $r_k$ orientable with respect to the homology theory over the field.

2) Each $A_k$ is a proper compact polyhedron subspace of $X_k$.
}\vspace{2mm}

{\it Proof}\, We have the following commutative diagram of exact sequences
\[\begin{array}{cccccccccc}
\cdots\longrightarrow \hspace{-1.5mm}&\hspace{-1.5mm} {\scriptstyle H_n(A_k)} \hspace{-1.5mm}&\hspace{-1.5mm} \stackrel{i_k}{\longrightarrow} \hspace{-1.5mm}&\hspace{-1.5mm} {\scriptstyle H_n(X_k)}
&\hspace{-1.5mm} \stackrel{j_k}{\longrightarrow} \hspace{-1.5mm}&\hspace{-1.5mm} {\scriptstyle H_n(X_k,A_k)} \hspace{-1.5mm}&\hspace{-1.5mm}\stackrel{\partial_k}{\longrightarrow} \hspace{-1.5mm}&\hspace{-1.5mm}
{\scriptstyle H_{n-1}(A_k)}\hspace{-1.5mm}&\hspace{-1.5mm} \longrightarrow \cdots\,\vspace{1mm}\\
&^{\alpha_k}\downarrow\quad&&^{\gamma_k}\downarrow\quad&&^{\beta_k}\downarrow\quad&&^{\alpha_k}\downarrow\quad\quad&&\\
\cdots\longrightarrow \hspace{-1.5mm}&\hspace{-1.5mm} {\scriptstyle H^{r_k-n}(X_k,A_k^c)} \hspace{-1.5mm}&\hspace{-1.5mm} \stackrel{q_k^\circ}{\longrightarrow} \hspace{-1.5mm}&\hspace{-1.5mm}{\scriptstyle H^{r_k-n}(X_k)}
\hspace{-1.5mm}&\hspace{-1.5mm} \stackrel{p_k^\circ}{\longrightarrow} \hspace{-1.5mm}&\hspace{-1.5mm}
{\scriptstyle H^{r_k-n}(A_k^c)} \hspace{-1.5mm}&\hspace{-1.5mm}\stackrel{\partial_k^\circ}{\longrightarrow} \hspace{-1.5mm}&\hspace{-1.5mm}
{\scriptstyle H^{r_k-n+1}(X_k,A^c_k)} \hspace{-1.5mm}&\hspace{-1.5mm}\longrightarrow\cdots,
  \end{array}\]
where $\alpha_k,\beta_k$ are the Alexander duality isomorphisms
and $\gamma_k$ is the Poncar\'{e} duality isomorphism.
So we have the following group isomorphisms
$$\begin{array}{rl}
(\partial_k^\circ)^{-1}\alpha_k\colon&\,\,\,{\rm ker}\,i_k\,\,\,\stackrel{\cong}
{\longrightarrow}\,\,\Sigma{\rm coker}\,p_k^\circ,\vspace{1mm}\\
\gamma_k\colon&\,\,\,{\rm im}\,i_k\,\,\,\,\stackrel{\cong}{\longrightarrow}
\,\,\,{\rm ker}\,p_k^\circ,\vspace{1mm}\\
p_k^\circ\gamma_k\colon&{\rm coker}\,i_k\stackrel{\cong}{\longrightarrow}
\,\,\,\,{\rm im}\,p_k^\circ.
\end{array}$$
Define $\theta_k\colon H_*^\XX(X_k,A_k)\to H^*_\XX(X_k,A_k^c)$
to be the direct sum of the above three isomorphisms.
Then $\theta_1{\otimes}{\cdots}{\otimes}\theta_m=\oplus_{(\sigma\!,\omega)\in\XX_m}\,\gamma_{\sigma\!,\omega}$.
The degree correspondence of $\theta_k$ is $H_*^{\emptyset,\emptyset}\to H^{r_k-*}_{\emptyset,\emptyset}$,
$H_*^{\{1\},\emptyset}\to H^{r_k-*}_{\{1\},\emptyset}$,
$H_*^{\emptyset,\{1\}}\to H^{r_k-*-1}_{\emptyset,\{1\}}$.

Now we construct similar homomorphism $\w\theta_k$ for the proof of Theorem~5.6. We have group isomorphisms
$$\begin{array}{rl}
\alpha_k\colon&\,\,\,{\rm ker}\,i_k\,\,\,\stackrel{\cong}
{\longrightarrow}\,\,{\rm im}\,\partial_k^\circ\subset H^*(X_k,A_k^c),\vspace{1mm}\\
\gamma_k\colon&\,\,\,{\rm im}\,i_k\,\,\,\,\stackrel{\cong}{\longrightarrow}
\,\,\,{\rm coim}\,q_k^\circ\subset H^*(X_k,A_k^c),\vspace{1mm}\\
p_k^\circ\gamma_k\colon&{\rm coker}\,i_k\stackrel{\cong}{\longrightarrow}
\,\,\,\,{\rm coim}\,p_k^\circ\subset H^*(X_k).
\end{array}$$
Define $\w\theta_k\colon H_*^\XX(X_k,A_k)\to\w H^*_\XX(X_k,A_k^c)$
to be the direct sum of the above three isomorphisms with $\w H^*_\XX(X_k,A_k)$ the direct sum
of the right side groups.
Then $\w\theta_1{\otimes}{\cdots}{\otimes}\w\theta_m=\oplus_{(\sigma\!,\omega)\in\XX_m}\,\w\gamma_{\sigma\!,\omega}$.
Since
$\w H^*_{\emptyset,\emptyset}\cong H^{*}_{\emptyset,\emptyset}$,
$\w H^*_{\{1\},\emptyset}\cong H^*_{\{1\},\emptyset}$,
$\w H^{*-1}_{\emptyset,\{1\}}\cong H^*_{\emptyset,\{1\}}$,
we have $\tilde\gamma_{\sigma\!,\,\omega}=\gamma_{\sigma\!,\,\omega}$.
\hfill $\Box$\vspace{3mm}

{\bf Theorem 5.2} {\it Let $K$ and $K^\circ$ be the dual of each other relative to $[m]$.
Then  for all $(\sigma\!,\omega)\in\XX_m$, $\omega\neq\emptyset$,
there are duality isomorphisms
$$\gamma_{K,\sigma\!,\,\omega}\colon H_*^{\sigma\!,\,\omega}(K)=\w H_{*-1}(K_{\sigma\!,\omega})\to
H^{|\omega|-*-1}_{\tilde\sigma\!,\,\omega}(K^\circ)=\w H^{|\omega|-*-2}((K^\circ)_{\tilde\sigma\!,\omega}),$$
$$\gamma^\circ_{K,\sigma\!,\,\omega}\colon H^*_{\sigma\!,\,\omega}(K)=\w H^{*-1}(K_{\sigma\!,\omega})
\to H_{|\omega|-*-1}^{\tilde\sigma\!,\,\omega}(K^\circ)=\w H_{|\omega|-*-2}((K^\circ)_{\tilde\sigma\!,\omega}),$$
where $\w\sigma=[m]{\setminus}(\sigma{\cup}\omega)$, $|\omega|$ is the cardinality of $\omega$.
}\vspace{2mm}

{\it Proof}\, Let $(C_*(\Delta\!^{\omega},K_{\sigma\!,\,\omega}),d)$ be the relative simplicial chain complex.
Since $\w H_*(\Delta\!^{\omega})=0$, we have a boundary isomorphism\vspace{1mm}\\
\hspace*{30mm}$\partial\colon H_*(\Delta\!^\omega,K_{\sigma\!,\,\omega})\stackrel{\cong}{\longrightarrow}
\w H_{*-1}(K_{\sigma\!,\,\omega})=H_*^{\sigma\!,\,\omega}(K)$.\vspace{1mm}\\
$C_*(\Delta\!^{\omega},K_{\sigma\!,\,\omega})$ has a set of generators consisting of all non-simplices of $K_{\sigma\!,\,\omega}$, i.e., $K^c_{\sigma\!,\,\omega}=\{\eta\subset\omega\,|\,\eta\not\in K_{\sigma\!,\,\omega}\}$ is a set of generators of $C_*(\Delta\!^{\omega},K_{\sigma\!,\,\omega})$.
So we may denote $(C_*(\Delta\!^\omega,K_{\sigma\!,\,\omega}),d)$ by $(C_*(K^c_{\sigma\!,\,\omega}),d)$, where
$\eta\in K^c_{\sigma\!,\,\omega}$ has degree $|\eta|{-}1$ with $|\eta|$ the cardinality of $\eta$.
The correspondence $\eta\to \omega{\setminus}\eta$ for all $\eta\in K^c_{\sigma\!,\,\omega}$ induces a dual complex isomorphism\vspace{1mm}\\
\hspace*{32mm}$\psi\colon(C_*(K^c_{\sigma\!,\,\omega}),d)\to
(\w C^{*}((K_{\sigma\!,\,\omega})^\circ),\delta)$.\vspace{1mm}\\
Since $(K_{\sigma\!,\,\omega})^\circ=(K^\circ)_{\tilde\sigma\!,\,\omega}$, we have induced homology group isomorphism
$\bar\psi\colon H_*(\Delta\!^\omega,K_{\sigma\!,\,\omega})\to H^{|\omega|-*-1}_{\tilde\sigma\!,\,\omega}(K^\circ)$.
Define $\gamma_{K,\sigma\!,\,\omega}=\bar\psi\partial^{-1}$.
\hfill $\Box$\vspace{3mm}

Notice that when $\omega=\emptyset$, for $\sigma\in K$, $\widetilde\sigma=[m]{\setminus}\sigma\notin K^\circ$
and there is no isomorphism from $H_*^{\sigma,\emptyset}(K)={\Bbb Z}$ to $H^*_{\tilde\sigma,\emptyset}(K^\circ)=0$;
for $\sigma\notin K$, $\widetilde\sigma=[m]{\setminus}\sigma\in K^\circ$
and there is no isomorphism from $H_*^{\sigma,\emptyset}(K)=0$ to $H^*_{\tilde\sigma,\emptyset}(K^\circ)={\Bbb Z}$.
\vspace{3mm}

{\bf Example 5.3} For the ${\cal S}(K;L_1,{\cdots},L_m)$ and index sets $\sigma\!,\,\omega,\hat\sigma,\hat\omega,\sigma_k,\omega_k$ in Example 3.12,
$\gamma_{{\cal S}(K;L_1,{\cdots},L_m),\,\sigma,\,\omega}
=\gamma_{K,\hat\sigma\!,\,\hat\omega}\otimes(\otimes_{\omega_k\neq\emptyset}\,\gamma_{L_k,\sigma_k,\omega_k})$.
\vspace{3mm}

{\bf Definition 5.4} For homology split $M={\cal Z}(K;\underline{X},\underline{A})$, let
$i\colon H_*(M)\to H_*(\w X)$ and $i^\circ\colon H^*(\w X)\to H^*(M)$
be the singular (co)homology homomorphism induced by the inclusion map from $M$ to
$\w X=X_1{\times}{\cdots}{\times}X_m$.
From the long exact exact sequences
$$\begin{array}{l}
\scriptstyle{\cdots\,\,\,\,\longrightarrow\,\,\,\, H_n(M)\,\,\,\,\stackrel{i}{\longrightarrow}\,\,\,\,H_n(\w X)
\,\,\,\,\stackrel{j}{\longrightarrow}\,\,\,\, H_n(\w X,M)\,\,\,\,\stackrel{\partial}{\longrightarrow}\,\,\,\, H_{n-1}(M)
\,\,\,\,\longrightarrow\,\,\,\,\cdots}\\\\
\scriptstyle{\cdots\,\,\,\,\longrightarrow\,\,\,\, H^{n-1}(M)\,\,\,\,\stackrel{\partial^\circ}{\longrightarrow}\,\,\,\, H^n(\w X,M)
\,\,\,\,\stackrel{j^\circ}{\longrightarrow}\,\,\,\, H^n(\w X)\,\,\,\,\stackrel{i^\circ}{\longrightarrow}\,\,\,\, H^{n}(M)
\,\,\,\, \longrightarrow\,\,\,\,\cdots}
  \end{array}$$
we define
$$\hat H_*(M)\!=\!{\rm coim}\,i,\,\overline H_*(M)\!=\!{\rm ker}\,i,\,
\hat H_*(\w X{,}M)\!=\!{\rm im}\,j,\,\overline H_*(\w X{,}M)\!=\!{\rm coker}\,j,$$
$$\hat H^*(M)\!=\!{\rm im}\,i^\circ,\,\overline H^*(M)\!=\!{\rm coker}\,i^\circ,\,
\hat H^*(\w X{,}M)\!=\!{\rm coim}\,j^\circ,\,H^*(\w X{,}M)\!=\!{\rm ker}\,j^\circ.\vspace{3mm}$$

{\bf Theorem 5.5} {\it For a homology split space $M={\cal Z}(K;\underline{X},\underline{A})$,
we have the following group decompositions
$$\begin{array}{cc}
H_*(M)=\hat H_*(M){\oplus}\overline H_*(M),&
H_*(\w X,M)=\hat H_*(\w X,M){\oplus}\overline H_*(\w X,M),\vspace{2mm}\\
H^*(M)=\hat H^*(M){\oplus}\overline H^*(M),&
H^*(\w X,M)=\hat H^*(\w X,M){\oplus}\overline H^*(\w X,M)
\end{array}$$
and direct sum group decompositions
$$\begin{array}{c}
\overline H_{*+1}(\w X,M)\cong \overline H_*(M)
\cong \oplus_{(\sigma\!,\,\omega)\in\overline\XX_m}\,
H_*^{\sigma\!,\,\omega}(K)\otimes H_*^{\sigma\!,\,\omega}(\underline{X},\underline{A}),\vspace{2mm}\\
\overline H^{*+1}(\w X,M)\cong \overline H^*(M)
\cong \oplus_{(\sigma\!,\,\omega)\in\overline\XX_m}\,
H^*_{\sigma\!,\,\omega}(K)\otimes H^*_{\sigma\!,\,\omega}(\underline{X},\underline{A}),\vspace{2mm}\\
\hat H_*(M)\cong\oplus_{\sigma\in K} H_*^{\sigma,\emptyset}(\underline{X},\underline{A}),\,\,
\hat H_*(\w X,M)\cong\oplus_{\sigma\notin K} H_*^{\sigma,\emptyset}(\underline{X},\underline{A}),\vspace{2mm}\\
\hat H^*(M)\cong\oplus_{\sigma\in K} H^*_{\sigma,\emptyset}(\underline{X},\underline{A}),
\,\,\hat H^*(\w X,M)\cong\oplus_{\sigma\notin K} H^*_{\sigma,\emptyset}(\underline{X},\underline{A}),
\end{array}$$
where $\overline\XX_m=\{(\sigma\!,\omega)\in\XX_m\,|\,\omega\neq\emptyset\}$.

The conclusion holds for all polyhedral product spaces if the (co)homology group is taken over a field.
}\vspace{3mm}

{\it Proof}\, By definition, $i=\oplus_{(\sigma\!,\,\omega)\in\XX_m}\,i_{\sigma\!,\,\omega}$ with\vspace{2mm}\\
\hspace*{15mm}$i_{\sigma\!,\,\omega}\colon H_*^{\sigma\!,\,\omega}(K){\otimes}H_*^{\sigma\!,\,\omega}(\underline{X},\underline{A})
\stackrel{i{\otimes}1}{-\!\!\!\longrightarrow}
H_*^{\sigma\!,\,\omega}(\Delta\!^{[m]}){\otimes}H_*^{\sigma\!,\,\omega}(\underline{X},\underline{A})$,\vspace{1mm}\\
where $i$ is induced by inclusion and $1$ is the identity.
So
$$\hat H_*(M)=\oplus_{\sigma\in K}
H_*^{\sigma,\emptyset}(K){\otimes}H_*^{\sigma,\emptyset}(\underline{X},\underline{A})
\cong\oplus_{\sigma\in K} H_*^{\sigma,\emptyset}(\underline{X},\underline{A})$$
$$\overline H_*(M)
=\oplus_{(\sigma\!,\,\omega)\in\overline\XX_m}\,
H_*^{\sigma\!,\,\omega}(K)\otimes H_*^{\sigma\!,\,\omega}(\underline{X},\underline{A}).$$

The relative group case is similar.
\hfill $\Box$\vspace{3mm}

{\bf Theorem 5.6} {\it For $M={\cal Z}(K;\underline{X},\underline{A})$ such that $(\underline{X},\underline{A})$
satisfies the condition of Theorem 5.1,
the Alexander duality isomorphisms
$$\alpha\colon H_*(M)\to H^{r-*}(\w X,M^c),\,\quad
\alpha^\circ\colon H^*(M)\to H_{r-*}(\w X,M^c)$$
have direct sum decomposition $\,\alpha=\hat\alpha\oplus\overline\alpha$,\,\,$\alpha^\circ=\hat\alpha^\circ\oplus\overline\alpha^\circ$,
where
$$\hat\alpha\colon \hat H_*(M)\to \hat H^{r-*}(\w X,M^c),\,\quad
\overline\alpha\colon \overline H_*(M)\to \overline H^{r-*}(\w X,M^c)\cong\overline H^{r-*-1}(M^c),$$
$$\hat\alpha^\circ\colon \hat H^*(M)\to \hat H_{r-*}(\w X,M^c),\,\quad
\overline\alpha^\circ\colon \overline H^*(M)\to \overline H_{r-*}(\w X,M^c)\cong\overline H_{r-*-1}(M^c)\,$$
are as follows. Identify all the above groups with the direct sum groups in Theorem~5.5.
Then
$$\hat\alpha\,\,=\oplus_{\sigma\in K}\,\gamma_{\sigma,\emptyset},\,\quad
\overline\alpha\,\,=\oplus_{(\sigma\!,\omega)\in\overline\XX_m}\,\,\gamma_{K,\sigma\!,\,\omega}{\otimes}\gamma_{\sigma\!,\,\omega},$$
$$\hat\alpha^\circ=\oplus_{\sigma\in K}\,\gamma_{\sigma,\emptyset}^\circ,\,\quad
\overline\alpha^\circ\,\,=\oplus_{(\sigma\!,\omega)\in\overline\XX_m}\,\,\gamma_{K,\sigma\!,\,\omega}^\circ{\otimes}\gamma_{\sigma\!,\,\omega}^\circ,$$
where $\gamma_{-}$, $\gamma_{-}^\circ$ are as in Theorem 5.1 and Theorem 5.2.
}\vspace{2mm}

{\it Proof}\, Denote by $\alpha=\alpha_M$, $\hat\alpha=\hat\alpha_M$,
$\overline\alpha=\overline\alpha_M$.
Then for $M={\cal Z}(K;\underline{X},\underline{A})$ and $N={\cal Z}(L;\underline{X},\underline{A})$,
we have the following commutative diagrams of exact sequences
\[\,\,\,\,\begin{array}{ccccccc}
{\scriptstyle\cdots\,\,\longrightarrow}\!\!&\!\!{\scriptstyle H_k(M{\cap}N)}\!\!&\!\!{\scriptstyle\longrightarrow}
\!\!&\!\!{\scriptstyle H_k(M){\oplus}H_k(N)}
\!\!&\!\!{\scriptstyle\longrightarrow}\!\!&\!\!{\scriptstyle H_k(M{\cup}N)}\!\!&\!\!{\scriptstyle\longrightarrow\,\,\cdots}\vspace{1mm}\\
\!\!&\!\!^{\alpha_{M\cap N}}\downarrow\quad\quad\!\!&\!\!\!\!&\!\!^{\alpha_M\oplus\alpha_N}\downarrow\quad\quad\quad\!\!&\!\!\!\!&\!\!^{\alpha_{M\cup N}}\downarrow\quad\quad\!\!&\!\!\\
{\scriptstyle\cdots\,\,\longrightarrow}\!\!&\!\!{\scriptstyle H^{r-k}(\w X,(M{\cap}N)^c)}\!\!&\!\!{\scriptstyle\longrightarrow}\!\!&\!\!
{\scriptstyle H^{r-k}(\w X,M^c){\oplus}H^{r-k}(\w X,N^c)}
\!\!&\!\!{\scriptstyle\longrightarrow}\!\!&\!\!{\scriptstyle H^{r-k}(\w X,(M{\cup}N)^c)}\!\!&\!\!{\scriptstyle\longrightarrow\,\,\cdots}
  \end{array}\quad\,\,\,\,(1)\vspace{1mm}\]
\[\begin{array}{ccccccc}
{\scriptstyle0\quad\longrightarrow}\!\!&\!\!{\scriptstyle \hat H_k(M{\cap}N)}\!\!&\!\!{\scriptstyle\longrightarrow}\!\!&\!\!
{\scriptstyle\hat H_k(M){\oplus}\hat H_k(N)}
\!\!&\!\!{\scriptstyle\longrightarrow}\!\!&\!\!
{\scriptstyle\hat H_k(M{\cup}N;\underline{X},\underline{A})}\!\!&\!\!{\scriptstyle\longrightarrow\quad 0}
\vspace{1mm}\\
\!\!&\!\!^{\hat\alpha_{M\cap N}}\downarrow\quad\quad\!\!&\!\!\!\!&\!\!^{\hat\alpha_M\oplus\hat\alpha_N}\downarrow\quad\quad\quad\!\!&\!\!\!\!&\!\!^{\hat\alpha_{M\cup N}}\downarrow\quad\quad\!\!&\!\!\\
{\scriptstyle0\quad\longrightarrow}\!\!&\!\!{\scriptstyle\hat H^{r-k}(\w X,(M{\cap}N)^c)}\!\!&\!\!{\scriptstyle\longrightarrow}\!\!&\!\!
{\scriptstyle\hat H^{r-k}(\w X,M^c){\oplus}\hat H^{r-k}(\w X,N^c)}
\!\!&\!\!{\scriptstyle\longrightarrow}\!\!&\!\!{\scriptstyle\hat H^{r-k}(\w X,(M{\cup}N)^c)}\!\!&\!\!{\scriptstyle\longrightarrow\quad 0}
  \end{array}\quad(2)\vspace{1mm}\]
For $(\sigma\!,\,\omega)\in\overline\XX_m$,
$A=H_l^{\sigma\!,\,\omega}(\underline{X},\underline{A})$,
$B=H^{r-|\omega|-l}_{\sigma\!,\,\omega}(\underline{X},\underline{A}^c)$,
$\gamma_1=\gamma_{K{\cap}L,\sigma\!,\omega}$, $\gamma_2=\gamma_{K,\sigma\!,\omega}{\oplus}\gamma_{L,\sigma\!,\omega}$,
$\gamma_3=\gamma_{K\cup L,\sigma\!,\omega}$, we have the commutative diagram
\[\begin{array}{ccccccc}
{\scriptstyle\cdots\,\,\longrightarrow}\!\!&\!\!{\scriptstyle H_k^{\sigma\!,\omega}(K{\cap}L){\otimes}A}\!\!&\!\!{\scriptstyle\longrightarrow}\!\!&\!\!
{\scriptstyle (H_k^{\sigma\!,\omega}(K){\oplus} H_k^{\sigma\!,\omega}(L)){\otimes}A}\!\!&\!\!
{\scriptstyle\longrightarrow}\!\!&\!\!{\scriptstyle H_k^{\sigma\!,\omega}(K{\cup}L){\otimes}A}\!\!&\!\!{\scriptstyle\longrightarrow\,\,\cdots}\vspace{1mm}\\
\!\!&\!\!^{\gamma_1\otimes\gamma_{\sigma\!,\,\omega}}
\downarrow\quad\quad\!\!&\!\!\!\!&\!\!^{\gamma_2\otimes\gamma_{\sigma\!,\,\omega}}
\downarrow\quad\quad\quad\!\!&\!\!\!\!&\!\!
^{\gamma_{3}{\otimes}\gamma_{\sigma\!,\,\omega}}\downarrow\quad\quad\!\!&\!\!\\
{\scriptstyle\cdots\,\,\longrightarrow}\!\!&\!\!
{\scriptstyle H^{|\omega|-k-1}_{\tilde\sigma\!,\omega}((K\cap L)^\circ){\otimes}B}\!\!&\!\!
{\scriptstyle\longrightarrow}\!\!&\!\!{\scriptstyle (H^{|\omega|-k-1}_{\tilde\sigma\!,\omega}(K^\circ)
{\oplus}H^{|\omega|-k-1}_{\tilde\sigma\!,\omega}(L^\circ)){\otimes}B}\!\!&\!\!{\scriptstyle\longrightarrow}\!\!&\!\!
{\scriptstyle H^{|\omega|-k-1}_{\tilde\sigma\!,\omega}((K\cup L)^\circ){\otimes}B}\!\!&\!\!
{\scriptstyle\longrightarrow\,\,\cdots}
  \end{array}\]
The direct sum of all the above diagrams is the following diagram.
\[\begin{array}{ccccccc}
{\scriptstyle\cdots\,\,\longrightarrow}\!\!&\!\!{\scriptstyle\overline H_k(M{\cap}N)}\!\!&\!\!{\scriptstyle\longrightarrow}\!\!&\!\!
{\scriptstyle\overline H_k(M){\oplus}\overline H_k(N)}\!\!&\!\!{\scriptstyle\longrightarrow}\!\!&\!\!
{\scriptstyle\overline H_k(M{\cup}N)}\!\!&\!\!{\scriptstyle\longrightarrow\,\,\cdots}\vspace{1mm}\\
\!\!&\!\!^{\overline\alpha_{M\cap N}}\downarrow\quad\quad\!\!&\!\!\!\!&\!\!^{\overline\alpha_M\oplus\overline\alpha_N}\downarrow\quad\quad\quad\!\!&\!\!\!\!&\!\!^{\overline\alpha_{M\cup N}}\downarrow\quad\quad\!\!&\!\!\\
{\scriptstyle\cdots\,\,\longrightarrow}\!\!&\!\!{\scriptstyle\overline H^{r-k}(\w X,(M{\cap}N)^c)}\!\!&\!\!
{\scriptstyle\longrightarrow}\!\!&\!\!{\scriptstyle\overline H^{r-k}(\w X,M^c){\oplus}\overline H^{r-k}(\w X,N^c)}
\!\!&\!\!{\scriptstyle\longrightarrow}\!\!&\!\!{\scriptstyle\overline H^{r-k}(\w X,(M{\cup}N)^c)}\!\!&\!\!
{\scriptstyle\longrightarrow\,\,\cdots}
  \end{array}\quad(3)\]

(1), (2) and (3) imply that
if the theorem holds for $M$ and $N$ and $M{\cap}N$, then it holds for $M{\cup}N$.
So by induction on the number of maximal simplices of $K$,
we only need prove the theorem for the special case that $K$ has only one maximal simplex.

Now we prove the theorem for $M={\cal Z}(\Delta\!^S;\underline{X},\underline{A})$ with $S\subset[m]$.
Then
$$M=Y_1{\times}{\cdots}{\times}Y_m,\quad Y_k=\left\{\begin{array}{cc}
X_k & {\rm if}\,\,k\in S, \\
A_k & {\rm if}\,\,k\notin S.
\end{array}\right.$$
So $(\w X,M^c)=(X_1,Y_1^c){\times}{\cdots}{\times}(X_m,Y_m^c)$.

Let $\w H^*_{\XX}(X_k,A_k),\w\theta_k,\gamma_{\sigma\!,\,\omega},\w\gamma_{\sigma\!,\,\omega}$ be as in the proof of Theorem 5.1.
Then we have the following commutative diagram
\[\begin{array}{ccc}
{\scriptstyle H_*(M)}
&\stackrel{\alpha_M}{-\!\!\!-\!\!\!-\!\!\!-\!\!\!-\!\!\!\longrightarrow}&
{\scriptstyle H^{r-*}(\w X,M^c)}\,\vspace{1mm}\\
\|&&\|\\
{\scriptstyle H_*(Y_1){\otimes}{\cdots}{\otimes}H_*(Y_m)}
&\stackrel{\alpha_M}{-\!\!\!-\!\!\!-\!\!\!-\!\!\!-\!\!\!\longrightarrow}&
{\scriptstyle H^{r_1-*}(X_1,Y_1^c){\otimes}{\cdots}{\otimes}H^{r_m-*}(X_m,Y_m^c)}\,\vspace{1mm}\\
\|\wr&&\|\wr\\
\oplus_{\sigma\subset S,\,\omega{\cap}S=\emptyset}\,{\scriptstyle H_*^{\sigma\!,\,\omega}(\underline{X},\underline{A})}
&\stackrel{\oplus\,\tilde\gamma_{\sigma\!,\omega}}{-\!\!\!-\!\!\!-\!\!\!-\!\!\!-\!\!\!\longrightarrow}&
\oplus_{\tilde\sigma\subset S,\,\omega{\cap}S=\emptyset}\,{\scriptstyle
\w H^{r-*}_{\tilde \sigma\!,\,\omega}(\underline{X},\underline{A}^c)}\vspace{1mm}\\
\cap&&\cap\\
{\scriptstyle H_*^\XX(X_1,A_1){\otimes}{\cdots}{\otimes}H_*^\XX(X_m,A_m)}
&\stackrel{\tilde\theta_1\otimes\cdots\otimes\tilde\theta_m}
{-\!\!\!-\!\!\!-\!\!\!-\!\!\!-\!\!\!\longrightarrow}&
{\scriptstyle \w H^*_\XX(X_1,A_1^c){\otimes}{\cdots}{\otimes}\w H^*_\XX(X_m,A_m^c)},
\end{array}\quad(4)\]
where $\w H^*_\XX(X_1,A_1^c){\otimes}{\cdots}{\otimes}\w H^*_\XX(X_m,A_m^c)=
\oplus_{(\sigma\!,\,\omega)\in\XX_m}\w H^{*}_{\sigma\!,\,\omega}(\underline{X},\underline{A}^c)$.

For $\sigma\subset S$, $\omega{\cap}S=\emptyset$, $\omega\neq\emptyset$,
we have
\[\begin{array}{ccc}
{\scriptstyle H_0^{\sigma\!,\,\omega}(\Delta\!^S){\otimes}H_*^{\sigma\!,\,\omega}(\underline{X},\underline{A})}
&\stackrel{\gamma_{\Delta\!^S\!,\sigma\!,\,\omega}\otimes\gamma_{\sigma\!,\,\omega}}
{-\!\!\!-\!\!\!-\!\!\!-\!\!\!-\!\!\!\longrightarrow}&
{\scriptstyle H^{|\omega|-1}_{\tilde\sigma\!,\,\omega}((\Delta\!^S)^\circ){\otimes}
H^{r-|\omega|-*}_{\tilde\sigma\!,\,\omega}(\underline{X},\underline{A}^c)\quad(\subset\overline H^{r-*-1}(M^c))}
\,\vspace{1mm}\\
\|\wr&&\hspace*{-20mm}\|\wr\\
{\scriptstyle H_*^{\sigma\!,\,\omega}(\underline{X},\underline{A})}
&\stackrel{\w\gamma_{\sigma\!,\,\omega}}{-\!\!\!-\!\!\!-\!\!\!-\!\!\!-\!\!\!\longrightarrow}&
{\scriptstyle
\w H^{r-*}_{\tilde \sigma\!,\,\omega}(\underline{X},\underline{A}^c)\quad(\subset\overline H^{r-*}(\w X,M^c))}.
\end{array}\]
So with the identification of the theorem, the third row of (4) is the direct sum
$\overline\alpha_M=\oplus_{\sigma\subset S,\,\omega{\cap}S=\emptyset,\,\omega\neq\emptyset}\,
\gamma_{\Delta\!^S\!,\sigma\!,\,\omega}{\otimes}\gamma_{\sigma\!,\,\omega}$,
$\hat\alpha_M=\oplus_{\sigma\subset S}\,\gamma_{\sigma,\emptyset}$ ($\w\gamma_{\sigma,\emptyset}=\gamma_{\sigma,\emptyset}$) and
$\alpha_M=\hat\alpha_M{\oplus}\overline\alpha_M$.\hfill $\Box$\vspace{3mm}

{\bf Example 5.7}\, Regard $S^{r+1}$ as one-point compactification of $\Bbb R^{r+1}$.
Then for $q\leqslant r$, the standard  space pair $(S^{r+1},S^q)$ is given by\vspace{1mm}\\
\hspace*{2mm}$S^q=\{(x_1,{\cdots},x_{r+1})\in\Bbb R^{r+1}\subset S^{r+1}
\,|\,x^2_1{+}{\cdots}{+}x_{q+1}^2=1,\,x_i=0,\,\,{\rm if}\,\,i>q{+}1\}.$

Let $M={\cal Z}_K\Big(\!
\begin{array}{ccc}
\scriptstyle{r_1{+}1} \!&\!\scriptstyle{\cdots}\!&\!\scriptstyle{r_m{+}1}\\
\scriptstyle{ q_1}   \!&\!\scriptstyle{\cdots}\!&\! \scriptstyle{q_m}\end{array}\!\Big)
={\cal Z}(K;\underline{X},\underline{A})$ be the polyhedral product space
such that $(X_k,A_k)=(S^{r_k+1},S^{q_k})$.
Since  $S^{r-q}$ is a deformation retract of $S^{r+1}{\setminus}S^q$,
the complement space $M^c={\cal Z}(K^\circ;\underline{X},\underline{A}^c)$
is homotopy equivalent to ${\cal Z}_{K^\circ}\Big(
\begin{array}{ccc}
\scriptstyle{r_1+\,1\,} &\scriptstyle{\cdots}&\scriptstyle{r_m+\,1\,}\\
\scriptstyle{r_1{-}q_1}   &\scriptstyle{\cdots}& \scriptstyle{r_m{-}q_m}\end{array}\Big)$.\vspace{1mm}

Since all $H_*^{\sigma\!,\,\omega}(\underline{X},\underline{A})\cong {\Bbb Z}$,
we may identify $H_*^{\sigma\!,\,\omega}(K){\otimes}H_*^{\sigma\!,\,\omega}(\underline{X},\underline{A})$
with $\Sigma^t H_*^{\sigma\!,\,\omega}(K)$,
where $t=\Sigma_{k\in\sigma}(r_k{+}1)+\Sigma_{k\in\omega}\,q_k$.
For $\sigma\subset[m]$, let ${\Bbb Z}_\sigma$ be the free group generated by
$\sigma$ with degree $0$. Then
\[\hat H_*(M)=\oplus_{\sigma\in K}\,\Sigma^{\Sigma_{k\in\sigma}\,(r_k{+}1)}{\Bbb Z}_\sigma,\]
\[\overline H_*(M)=\oplus_{(\sigma\!,\,\omega)\in\overline \XX_m}\,
\Sigma^{\Sigma_{k\in\sigma}\,(r_k{+}1)+\Sigma_{k\in\omega}\,q_k}
H_*^{\sigma\!,\,\omega}(K).\]

Dually, the cohomology of the complement space $M^c$ is
\[\hat H^*(M^c)=\oplus_{\sigma\in K^\circ}\,\Sigma^{\Sigma_{k\in\sigma}\,(r_k{+}1)}{\Bbb Z}_\sigma,\]
\[\overline H^*(M^c)=\oplus_{(\sigma\!,\,\omega)\in\overline \XX_m}\,
\Sigma^{\Sigma_{k\in\sigma}\,(r_k{+}1)+\Sigma_{k\in\omega}\,(r_k-q_k)}
H^*_{\sigma\!,\,\omega}(K^\circ).\]

In this case, the direct sum of $\gamma_{K,\sigma\!,\omega}
\colon H_*^{\sigma\!,\,\omega}(K)\to H^{|\omega|-*-1}_{\tilde\sigma\!,\,\omega}(K^\circ)$
over all $(\sigma\!,\omega)\in\overline\XX_m$ (regardless of degree) is the isomorphism
$\overline H_*(M)\cong\overline H\,^{r-*-1}(M^c)$.

Specifically, ${\cal Z}(K;S^{2n+1},S^n)={\cal Z}_K\Big(\!
\begin{array}{ccc}
\scriptstyle{2n{+}1} \hspace{-1mm}&\hspace{-1mm}\scriptstyle{\cdots}\hspace{-1mm}&\hspace{-1mm}\scriptstyle{2n{+}1}\\
\scriptstyle{n}   \hspace{-1mm}&\hspace{-1mm}\scriptstyle{\cdots}\hspace{-1mm}&\hspace{-1mm} \scriptstyle{n}\end{array}\!\Big)$.
Then we have $$\overline H_*({\cal Z}(K;S^{2n+1},S^n))\cong
\overline H\,^{(2n+1)m-*-1}({\cal Z}(K^\circ;S^{2n+1},S^n)).$$

{\bf Acknowledgement} The author expresses his deepest gratitude to the referee for giving so many helpful
advices.

\end{document}